\documentclass[12pt, a4paper]{amsart}

\usepackage[english]{babel}
\usepackage{amsfonts}
\usepackage{amsmath}
\usepackage{amssymb, latexsym}
\usepackage{graphicx}
\usepackage[usenames]{color}

\usepackage{fullpage}
\usepackage{amsmath}

\usepackage{mathrsfs}

\newtheorem{theorem}{Theorem}[section]
\newtheorem{proposition}[theorem]{Proposition}
\newtheorem{corollary}[theorem]{Corollary}
\newtheorem{lemma}[theorem]{Lemma}
\newtheorem{remark}[theorem]{Remark}
\newtheorem{definition}[theorem]{Definition}

\numberwithin{equation}{section} 

\def\cala{{\mathcal{A}}}
\def\cale{{\mathcal{E}}}
\def\calf{{\mathcal{F}}}
\def\caln{\mathcal{N}}
\def\calo{\mathcal{O}}
\def\calh{\mathcal{H}}
\def\calv{\mathcal{U}}
\def\calv{\mathcal{V}}

\def\calb{{\mathcal{B}}}

\def\calu{{\mathcal{U}}}
\def\bbc{{\mathbb{C}}}
\def\bbe{{\mathbb{E}}}
\def\bbr{{\mathbb{R}}}
\def\bbp{{\mathbb{P}}}
\def\bbq{{\mathbb{Q}}}
\def\bbn{{\mathbb{N}}}
\def\bbm{{\mathbb{M}}}

\def\en t{{{\rm Z}\mkern-5.5mu{\rm Z}}}

\def\<{\left<}
\def\>{\right>}
\def\({\left(}
\def\){\right)}

\def\D{\Delta}
\def\9{{\infty}}
\def\barr{\begin{array}}
\def\earr{\end{array}}

\def\wt{\widetilde}

\def\ol{\overline}

\def\vf{{\varphi}}
\def\lbb{{\lambda}}
\def\g{{\gamma}}
\def\a{{\alpha}}

\def\D{{\Delta}}

\def\3{\subset }
\def\na{{\nabla}}

\def\ve{{\varepsilon}}
\def\p{{\partial}}
\def\l<{{\langle}}
\def\r>{{\rangle}}

\begin{document}

\title{Optimal control of nonlinear stochastic
differential equations on Hilbert spaces}

\author{Viorel Barbu}
\address[V. Barbu]{Octav Mayer Institute of
Mathematics (Romanian Academy) and Al.I. Cuza University,
700506, Ia\c si, Romania. }
\author{Michael R\"ockner}
\address[M. R\"ockner]{Fakult\"at f\"ur
Mathematik, Universit\"at Bielefeld,  D-33501 Bielefeld, Germany,
and Academy of Mathematics and Systems Science,
CAS, Beijing, China.}
\author{Deng Zhang}
\address[D. Zhang]{Department of Mathematics,
Shanghai Jiao Tong University, 200240 Shanghai, China.}
\begin{abstract}
We here consider  optimal control problems
governed by nonlinear stochastic equations on a Hilbert space $H$
with nonconvex payoff,
which is rewritten as a deterministic optimal control problem governed by
a Kolmogorov equation in $H$.
We prove  the existence and first-order necessary condition
of closed loop optimal controls
for the above control problem.
The  strategy is based on
solving a deterministic bilinear optimal control problem
for the corresponding Kolmogorov equation on the space $L^2(H,\nu)$,
where $\nu$ is the related infinitesimally invariant measure
for the Kolmogorov operator.
\end{abstract}


\maketitle
\section{Introduction}

We are concerned with optimal control problems
connected with the informal
stochastic differential equation on a Hilbert space $H$
(with norm $|\cdot|_H$, inner product $\<\cdot, \cdot\>$) of type
\begin{align} \label{equa-X}
   &dX(t) = A(X(t))dt + Q^\frac 12 B u(X(t)) dt + Q^\frac 12 dW(t), \ \ t\in (0,T),\\
   &X(0) = x \in H.  \nonumber
\end{align}
Here,
the operator $A$ is defined by
\begin{align} \label{def-A}
  A:D(F) \mapsto (D(\wt A))^*, \
  \<A(x),h\>:=\l<x,\wt A h\r> + \<F(x), h\>
\end{align}
for any $x\in D(F)$ and any $h\in D(\wt A)$,
where $\wt{A}$ is a self-adjoint
m-dissipative linear operator in $H$ to be made precise later on,
and $F$ is a (possibly nonlinear) operator from  $D(F) \subseteq H$ to $H$.

The operator
$B$ is linear and bounded  on $L^\9(H; H,  \nu)$,
where $\nu$ is an infinitesimally invariant measure for
the corresponding Kolmogorov operator when $u\equiv 0$
(see Hypothesis (H1) (ii) below),
which actually serves as a substitute for Lebesque measure on $H$
that does not exist on infinite dimensional spaces.
The operator
$Q$ is a positive definite bounded self-adjoint linear operator on $H$,
satisfying that
$Qe_i = q_i e_i$, $q_i > 0$ for all but finitely many $i\geq 1$ and
for some orthonormal basis $\{e_i\}_{i \geq 1} \subseteq D(\wt A)$ of $H$,
and $W$ is a cylindrical Wiener process on $H$ defined on
a probability space $(\Omega, \mathscr{F}, \bbp)$ with normal filtration $(\mathscr{F}_t)$, $t\geq 0$.

The term $u$ is an input controller applied to the stochastic system
and is taken in the admissible set
\begin{align*}
   \calu_{ad} = \{ u: H \rightarrow H;\ u\ is\ \nu-measurable,\ |u(x)|_H \leq \rho,\ \forall x\in H\},
\end{align*}
where $\rho \in (0,\9)$ is fixed.

Equation \eqref{equa-X} is mainly motivated by a number of stochastic partial differential equations,
including singular stochastic equations  (\cite{DRW09, DR02}),
gradient systems,
stochastic reaction-diffusion equations (\cite{DDG02,DL14})
and
stochastic porous media equations (\cite{BBDR06}, see also \cite{BDR16}).

In the present work,
we are interested in the optimal feedback control problem
for \eqref{equa-X}, i.e.,
find a controller $u^*\in \calu_{ad}$ such that
\begin{enumerate}
  \item[$(P_0)$]
  \begin{align*}
  {\rm Min} \bigg\{&\bbe \int_0^T \int_H g(X^u(t,x)) \nu(dx) dt;\ u\in \calu_{ad},\
               X^u\ solves\ \eqref{equa-X} \bigg\},
  \end{align*}
\end{enumerate}
where $g$ is a given function in $L^2(H,\nu)$,
is attained at $u^*$.

It should be mentioned that
the main difficulty for the existence theory for the optimal control
problem $(P_0)$ is that
the cost functional
$\Phi(u) = \bbe \int_0^T \int_H g(X^u(t,x)) \nu(dx) dt$,
$u\in \calu_{ad}$,
is not weakly lower-semicontinuous
on $L^2(H; H,\nu)$,
if $A$ is nonlinear and $g$ is not convex.

Another delicate problem in infinite dimensional spaces is that,
even if \eqref{equa-X} has a unique strong solution
(in the probabilistic sense)
in the uncontrolled case where $u\equiv 0$,
it is in general not clear whether it still has strong solutions
under bounded perturbations.
See, e.g., \cite{DFPR13}-\cite{DFRV16} for the relevant work.

Here, the key idea
is to
rewrite the original Problem $(P_0)$
as a deterministic bilinear optimal control problem
governed by the Kolmogorov equation corresponding to \eqref{equa-X}.

More precisely,
we consider
the  Kolmogorov equation corresponding to \eqref{equa-X}, i.e.,
\begin{align} \label{equa-Nu}
    &\frac{d\vf}{dt} = N_2 \vf + \l<Q^\frac 12 B u, D\vf\r>,\ \ t> 0,  \\
    & \vf(0,x) = g(x),\ \ x\in H,    \nonumber
\end{align}
where
$u\in \mathcal{U}_{ad}$,
$N_2$ is the Kolmogorov operator related to \eqref{equa-X}
(see \eqref{N0} and Remark \ref{Rem-i-ii-iii} below),
and equation \eqref{equa-Nu} is taken in the space $L^2(H, \nu)$.

Heuristically,
via It\^o's formula,
one has that
the solution $\vf^u$ for \eqref{equa-Nu} is given by
\begin{align*}
   \vf^u(t,x)
   = \bbe g(X^u(t,x)),\ for\ dt\times \nu-a.e.\ (t,x)\in[0,T]\times H.
\end{align*}
This entails that
the original optimal control problem can be reformulated as follows:
find $u^* \in \calu_{ad}$ such that
\begin{enumerate}
  \item[$(P^*)$]
  \begin{align*}
     &{\rm Min} \bigg\{ \int_0^T \int_H \vf^u(t,x) \nu(dx)dt;\ u\in \mathcal{U}_{ad}, \
        and\ \vf^u\ is\ the\  solution\ to\ \eqref{equa-Nu} \bigg\}
  \end{align*}
\end{enumerate}
is attained at $u^*$.

This idea was recently applied in \cite{B19} by the first author
to the  stochastic reflection problem in finite dimensions.
The main advantage of Problem $(P^*)$ is that
it is a deterministic bilinear optimal control problem.
This feature
makes it possible to give a unified treatment of optimal control problems
for various stochastic equations on Hilbert spaces
through the corresponding Kolomogorov operators,
under unusually weak conditions of the nonlinearity and the objective functionals.
Actually,
the usual continuity or convexity conditions are not assumed here,
which can be viewed as a regularization effect of noise
on control problems through
the corresponding Kolmogorov operators.

As a matter of fact,
the optimal feedback controllers for Problem $(P)$
can be formally determined by solving an infinite dimensional second order
Hamilton-Jacobi equation
(see, e.g., \cite{CIL92,FGS17}).
However,
such an equation under quite restrictive conditions has only a viscosity solution
which is not sufficiently regular to provide an explicit
representation for the optimal controller.
We would also like to refer to \cite{FT02}
for the solvability of nonlinear Kolmogorov equations, 
including Hamilton-Jacobi-Bellman equations, 
and the applications to optimal feedback controls.  

Here,
for any objective functions $g$ in $D(N_2)$,
where $D(N_2)$ is the domain of the closure in $L^2(H,\nu)$ of the Kolmogorov operator $(N_0, D(N_0))$
defined in \eqref{N0} below,
we prove the existence of a closed-loop optimal control
for Problem $(P^*)$ under mild conditions on $F$ and $g$.

Moreover,
in the symmetric case
(i.e.,
$N^*_2  = N_2$ on $L^2(H,\nu)$,
where $N_2^*$ denotes the dual operator of $N_2$),
for more general objective functions $g\in L^2(H,\nu)$,
we obtain the existence as well as
first-order necessary condition of optimal feedback controllers
of Problem $(P^*)$.

Regarding the original control problem of the stochastic equation \eqref{equa-X},
it turns out that
the martingale problem serves as an appropriate concept of solutions
to stochastic equations on Hilbert spaces.
More precisely, we consider the problem
\begin{enumerate}
  \item[$(P)$]
  \begin{align*}
  {\rm Min} \bigg\{\int_0^T \int_H & \bbe_{\bbp_x} g(X^u(t)) \nu(dx) dt;\ u\in \calu_{ad},\
                      \bbp_x\circ (X^u)^{-1}\ solves\  \\
               & the\ martingale\ problem\ of\ \eqref{equa-X}\
                for\ \nu-a.e.\ x\in H \bigg\}.
  \end{align*}
\end{enumerate}
(See Definition \ref{Def-equacontr} below for the definition of the
martingale problem corresponding to \eqref{equa-X}.)

We prove that
the optimal controllers to Problem $(P^*)$ obtained above
actually coincide with those to the  Problem $(P)$,
as long as the related martingale problems are well posed.
In this sense,
the optimal controllers for Problem $(P^*)$ of Kolmogorov equations
can be viewed as generalized optimal controllers
for the  Problem $(P)$ of stochastic equations on Hilbert spaces.

Actually,
the solutions to the martingale problem for \eqref{equa-X}
suffice to define the objective functional in Problem $(P)$.
More importantly,
well-posedness for this type of martingale problem
holds in a quite general setting
(e.g. in the framework of (generalized) Dirichlet forms),
and it is also stable under bounded perturbations
and thus enables us to treat optimal control problems of stochastic differential equations
on Hilbert spaces,
of which the nonlinearity may be not continuous or
the operator $Q^\frac 12$ is not necessarily Hilbert-Schmidt
(see, e.g., \cite{DR02}).
For such equations,
it is known that strong  solutions (in the probabilistic sense)
do not exist in general.

As we shall see below,
the martingale problem is well posed for various stochastic equations on Hilbert spaces,
including  singular dissipative stochastic  equations,
stochastic reaction-diffusion equations
as well as stochastic porous media equations.
Moreover,
we also prove that
the well-posedness of martingale problems are implied by
the m-dissipativity of the corresponding Kolmogorov operators in certain situations,
by using the theory of (generalized) Dirichlet forms
(see Theorems \ref{Thm-H1iv} and \ref{Thm-H1'-H1-Sym} below).
The interplay between optimal control problems
and   (generalized) Dirichlet forms
would be of independent interest.

We would also like to mention that,
by the argument above,
the end point optimal control problem
\begin{align*}
  {\rm Min}\{  \int_H  & \bbe_{\bbp_x} g(X^u(T)) \nu(dx):\ u\in \calu_{ad},\
              \bbp_x\circ (X^u)^{-1}\ solves\ \\
              &  the\ martingale\ problem\ of\  \eqref{equa-X}\
                 for\ \nu-a.e.\ x\in H.\}
\end{align*}
can be also written as
\begin{align*}
  {\rm Min} \{ \int_H \vf^u(T,x) \nu(dx):\ u\in \calu_{ad},\ \vf^u\ solves\ \eqref{equa-Nu}\}.
\end{align*}

{\bf Notation}
For $k\in \mathbb{N}$,
by  $\mathcal{F} C_b^k(H)$ we denote
the set of $C_b^k$-cylindrical functions
$\vf(x) = \phi(\<x,e_1\>,\cdots, \<x,e_n\>)$
for some $n\in \mathbb{N}$ and
$\phi\in C_b^k (\bbr^n)$,
where $\{e_k: k\in \bbn\}$
is the eigenbasis of $Q$ introduced above.
Let $B_b(H)$ and $C_b(H)$ denote, respectively,
the bounded Borel-measurable
and bounded continuous functions from $H$ to $\bbr$,
and let $L(H)$
be the set of all bounded operators on $H$.
The symbols $D$ and $D^2$ denote the first and second Fr\'{e}chet derivatives, respectively.
We also use the notation $Id$ for the identity operator on $H$.

For any Borel probability measure $\nu$ on $H$,
${\rm supp}(\nu)$ denotes the topological support of $\nu$,
and $L^2(H, \nu)$ consists of $\nu$-measurable functions $\vf$ on $H$
such that $\int_H |\vf(x)|^2 \nu(dx)<\9$.
We use the notation $(\ ,\ )$ for the inner product in $L^2(H,\nu)$.
Similarly,
$L^2(H; H, \nu)$ denotes the space of
$H$-valued $L^2(\nu)$-integrable maps.

\section{Formulation of the main results} \label{Sec-Main}
To begin with,
let us first introduce
the Kolmogorov operator related to \eqref{equa-X},
which is formally given by,
\begin{align} \label{Nu0}
   N^u_0  \vf (x) :=&  \frac 12 {\rm Tr} [Q D^2 \vf] (x)
                   + \l<A(x), D \vf(x)\r> + \l<Bu(x), Q^\frac 12 D \vf(x)\r>,
\end{align}
for any $\vf  \in \mathcal{F}C_b^2(H)$.
In particular,
when $u\equiv 0$,
we set
\begin{align} \label{N0}
   N_0  \vf (x) :=&  \frac 12 {\rm Tr} [Q D^2 \vf] (x)
                   + \<A(x), D \vf(x)\>, \ \ \vf \in D(N_0):= \calf C_b^2(H).
\end{align}

Consider the following assumptions.
\begin{enumerate}
  \item[(H1)]
  There exists a Borel probability measure $\nu$ such that
  $F: D(F)\subseteq H \to H$ is $\nu$-measurable and the following properties hold:
  \begin{enumerate}
     \item[(i)] $\nu(D(A)) =1$
                and $\int_H (|F(x)|^2_H + |x|^2_H)\nu(dx)<\9$.
     \item[(ii)] $\nu$ is the infinitesimally invariant measure for $(N_0, D(N_0))$,
                 i.e.,
                 $$ \int_H N_0 \vf d\nu =0,\ \  \forall \vf \in \mathcal{F} C_b^2(H).$$
     \item[(iii)] $(N_0, \mathcal{F} C_b^2(H))$ is essentially m-dissipative on $L^2(H, \nu)$,
                 i.e., $(1-N_0) (\calf C_b^2(H))$ is dense in $L^2(H,\nu)$.
  \end{enumerate}
  \item[(H2)]  The operator $Q^\frac 12 B$ with domain
              $D(Q^\frac 12 B):= \calu_{ad}$ and
              defined by $(Q^\frac 12 B)$ $(u)(x) : = Q^\frac 12 (Bu(x))$, $x\in H$,
              is compact as an operator from $L^\9(H;H,\nu)$ to $L^2(H;H,\nu)$,
              i.e.,
              if $u_n , u\in D(Q^\frac 12 B)$, $n\in \mathbb{N}$,
              such that $u_n\to u$ weakly-star in $L^\9(H;H,\nu)$ as $n\to \9$,
              then $Q^\frac 12 Bu_n \to Q^\frac 12 Bu$ in $L^2(H;H,\nu)$.
  \item[(H3)] The  operator $Q^\frac 12 D$ with domain $\calf C_b^1(H)$ is closable from $L^2(H, \nu)$ to $L^2(H; H, \nu)$,
              and the embedding $ W^{1,2}(H,\nu)$ into $L^2(H, \nu)$ is compact.
\end{enumerate}
Here
$W^{1,2}(H,\nu)$
is the Sobolev space defined
as the completion of $\calf C^2_b(H)$ under the norm
$  \| \vf \|_{W^{1,2}(H,\nu)}
    =( \int_H (|\vf|^2+|Q^\frac 12 D\vf|_H^2) d\nu)^\frac 12$.
Note that
$W^{1,2}(H,\nu)$ is a subspace of $L^2(H,\nu)$
if and only if $(Q^\frac 12D, \calf C_b^1(H))$ is closable,
as an operator from $L^2(H,\nu)$ to $L^2(H; H,\nu)$.
In this case we denote its closure again by $Q^\frac 12D$
and by construction its domain is $W^{1,2}(H,\nu)$.

\begin{remark} \label{Rem-i-ii-iii}
As is well-known $(H1)$ $(ii)$ implies that
$(N_0, D(N_0))$ is dissipative on $L^2(H,\nu)$,
so by $(H1)$ $(iii)$
and the Lumer-Phillips Theorem
its closure $(N_2,D(N_2))$ generates a $C_0$-semigroup
$P_t^\nu = e^{tN_2}$, $t>0$,
of contractions on $L^2(H,\nu)$.
Furthermore, $D(N_0)$ is dense in $D(N_2)$
with respect to the graph norm given by $N_2$.
\end{remark}

\begin{remark}
The compactness of the embedding of $W^{1,2}(H,\nu)$  into $L^2(H,\nu)$
is equivalent to the compactness of the  semigroup $P^\nu_t$ for some (equivalently, all) $t>0$.
See, e.g., \cite[Theorem 1.2]{GW02},
\cite[Theorems 1.1 and 3.1]{W00}
and \cite[p.3250]{W17}.
The above compact embedding can be also deduced from the Logarithmic-Sobolev inequality,
see, e.g., \cite{DDG02}.
In particular,
Hypothesis $(H3)$ holds for the
Gaussian invariant measures
of the Ornstein-Uhlenbeck process
(see \cite{CG95}).
\end{remark}

Below we give one specific example satisfying Hypothesis $(H2)$.\\

{\bf Example}
Let $Q = Id$,
and let $f_j\in  L^1(H; H, \nu) $, $g_j \in L^\9(H; H, \nu)$,
$j\geq 1$,
be such that
\begin{align} \label{B}
     B(u)  = \sum\limits_{j=1}^\9 \int_H \langle u, f_j \rangle d \nu  g_j,\ \ \forall u\in \calu_{ad},
\end{align}
and
\begin{align} \label{fjvej-l1}
     C_B:= \sum\limits_{j=1}^\9 \|f_j\|_{L^1(H; H, \nu)} \|g_j\|_{L^\9(H; H, \nu)} <\9.
\end{align}
Then, $B$ satisfies $(H2)$.

In fact, let $u_n, u\in D(B)$, $n\in \mathbb{N}$,
be such that
$u_n\to u$ weakly-star in $L^\9(H; H,\nu)$ as $n\to \9$.
Then, for every $N\in \mathbb{N}$,
\begin{align*}
    \sum_{j=1}^N \int_H \<u-u_n, f_j\> d\nu g_j \to 0,
    \ \ in\ L^\9(H; H,\nu).
\end{align*}
So, let $\ve >0$.
Then, by \eqref{fjvej-l1} there exists $N\in \mathbb{N}$ such that
\begin{align*}
   \sum_{j= N+1}^\9 \|f_j\|_{L^1(H; H,\nu)} \| g_j\|_{L^\9(H; H,\nu)} <\ve.
\end{align*}
Hence for $\ve >0$,
\begin{align*}
   \limsup\limits_{n\to \9}  \|Bu-Bu_n\|_{L^\9(H; H,\nu)}
   \leq&  \limsup\limits_{n\to \9}
          \|\sum\limits_{j=1}^N \int_H \<u-u_n, f_j\> d\nu g_j \|_{L^\9(H;H,\nu)} \\
       &  + 2 \rho \sum_{j= N+1}^\9 \|f_j\|_{L^1(H; H,\nu)} \| g_j\|_{L^\9(H; H,\nu)} \\
   \leq & 2\rho \ve.
\end{align*}
So, we even have $Bu_n \to B u$ in $L^\9(H; H,\nu)$ as $n\to \9$.

\begin{remark}
A specific example is where
$B(\calu_{ad})$ is in a finite dimensional subspace of $L^\9(H; H, \nu)$.

Actually,
in this case, there exists
linear independent $g_1,\cdots, g_j\in L^\9(H; H, \nu)$
such that $\|g_i\|_{L^2(H;H,\nu)}=1$, $1\leq i\leq j$
and $\{B(\calu_{ad})\} \subseteq span\{g_1, \cdots, g_n\}$,
and so,
for any $u\in \calu_{ad}$,
$B(u) = \sum_{j=1}^n c_j g_j$
for some $c_j\in \bbr$, $1\leq j\leq n$.
Then, we take $\{\wt{g}_j\}_{j=1}^n \subseteq L^2(H; H,\nu)$ such that
$\<\wt{g}_j, g_k\>_{L^2(H; H, \nu)}= \delta_{jk}$, $1\leq j,k\leq n$.
This yields that
$c_j = \<B(u), \wt{g}_j\>_{L^2(H; H, \nu)} = \<u, B^* \wt{g}_j\>_{L^2(H; H, \nu)}$,
where $B^*$ is the dual operator of $B$ in $L^2(H;H,\nu)$.
This implies \eqref{B} with $f_j = B^* \wt{g}_j$.
\end{remark}

Under Hypothesis $(H1)$,
let $(N_2, D(N_2))$ be
the closure of $(N_0,$ $\mathcal{F} C_b^2(H))$ in $L^2(H, \nu)$.
Then, $\nu$ is an invariant measure for $P_t^\nu = e^{t N_2}$, $t>0$,
i.e.,
\begin{align}
   \int_H P^\nu_t f(x) \nu(dx) =  \int_H   f(x) \nu(dx),\ \ \forall t\geq 0, \ \forall f\in \calb_b(H).
\end{align}
(See, e.g., the proof of \cite[Corollary 5.3]{DR02}.)

The essential m-dissipativity of $(N^u_0, \mathcal{F} C_b^2(H))$
can be inherited from
the uncontrolled case where $u\equiv 0$,
more precisely,
from $(H1)$ $(iii)$.
This is the content of the following theorem to be proved in Section \ref{Sec-OP-KE} below.

\begin{theorem}  \label{Thm-Max-Nu}
Assume Hypothesis $(H1)$ to hold.
Then, we have the integration by parts formula
\begin{align}  \label{Integ-part}
   \int_H  \vf N_2 \vf d\nu
   = -\frac 12 \int_H  |\ol{Q^\frac 12 D} \vf|_H^2 d \nu,
   \ \ \forall \vf \in D(N_2),
\end{align}
where  $\ol{Q^\frac 12 D}$ is the continuous extension of the operator
\begin{align*}
   D(N_0) \ni \vf \mapsto Q^\frac 12 D \vf \in L^2(H;H,\nu)
\end{align*}
with respect to the $N_2$-graph norm on $D(N_0)$.

Moreover,
for each $u\in \calu_{ad}$,
the operator
\begin{align*}
   N_2^u: D(N_2) \mapsto L^2(H, \nu),
   \ \ N_2^u\vf(x) := N_2 \vf(x) + \l<B u(x), \ol{Q^\frac 12 D} \vf(x)\r>
\end{align*}
has $\mathcal{F} C_b^2(H)$ as a core
and generates a $C_0$-semigroup $ e^{tN_2^u} $ on $L^2(H, \nu)$.
Furthermore,
for some positive constant $C(T,\rho)>0$,
\begin{align} \label{bdd-eN2u}
    \sup\limits_{u\in \calu_{ad}}
    (\| e^{tN_2^u}  g\|_{C([0,T];L^2(H, \nu))}
    + (\int_0^T \int_H |\ol{Q^\frac 12 D}  e^{tN_2^u}  g|_H^2 d\nu dt)^\frac 12)
     \leq C(T,\rho) \|g\|_{L^2(H,\nu)}.
\end{align}
\end{theorem}

The first result of this paper is concerned with
the existence of optimal controllers for Problem $(P^*)$.
It will be proved in Section \ref{Sec-OP-KE} below.
\begin{theorem} \label{Thm-Contr}
({\it Optimal control of Kolmogorov equations: general case})

Assume that  Hypothesis $(H1)$ holds
and, in addition, that either $(H2)$ or $(H3)$ holds.
Then,
for any  $g\in D(N_2)$,
there exists at least one optimal control
to Problem $(P^*)$.

In particular,
in the case where $(H1)$ and $(H3)$ hold,
one may take $B=Id$.
\end{theorem}

\begin{remark}
We would like to mention that,
no continuity or convexity of $A$ and $g$
are assumed in Theorem \ref{Thm-Contr}
which, however,
are the usual conditions for optimal feedback controls
even in the finite dimensional case.
\end{remark}

Next,
we are concerned with the symmetric case
(i.e., $N_2^* = N_2$ on $L^2(H,\nu)$)
which
arises, in particular,
in various applications
to gradient systems (see, e.g., \cite{DFRV16,DL14} and
the end of Subsection \ref{Subsec-SingSDE} below).

In this case,
we are able to obtain optimal controllers for more general objective
functions $g\in \calh:= L^2(H, \nu)$.
Moreover,
we also obtain the first-order necessary condition of the optimal feedback controllers,
in terms of the solutions to Kolmogorov equations and adjoint backward equations.

One nice feature here is that
the corresponding Kolmogorov operators are
defined in the variational form from
$\calv:= W^{1,2}(H,\nu)$  to $\calv'$,
where $\calv'$ is the dual space of $\calv$ in the pairing $(\cdot, \cdot)$
with the pivot space $\calh := L^2(H,\nu)$.
(Note that,
$\calv \subset \calh \subset \calv'$
with dense and continuous embeddings.)
This fact enables us to analyze the Kolmogorov equations
and the adjoint backward equations
in the variational setting.

The following result generalizes Theorem \ref{Thm-Contr}
in the symmetric case for $g\in L^2(H,\nu)$
and is proved in Section \ref{Sec-OP-KE}.
\begin{theorem} \label{Thm-Contr-Symm}
({\it Optimal control of Kolmogorov equations: symmetric case})

Consider the symmetric case $N_2^* = N_2$ on $L^2(H,\nu)$.
Assume $(H1)$
and that $(Q^\frac 12 D, \mathcal{F}C_b^1(H))$ is closable from
$L^2(H,\nu)$ to $L^2(H;H,\nu)$.
In addition,
assume that either $(H2)$ or $(H3)$ holds.

Then,
for any objective function  $g\in L^2(H,\nu)$,
there exists an optimal control for Problem $(P^*)$
where $\vf^u$ solves the equation in the space $\calv'$.

In particular,
in the case where $(H3)$ holds,
we can take $B=Id$.
\end{theorem}

In order to identify the optimal feedback controllers,
we (in the symmetric case)
introduce
the adjoint-backward equation corresponding to the Kolmogorov equation \eqref{equa-Nu}
\begin{align} \label{equa-p}
   & \frac{dp}{dt} = -N_2 p  - G^{u} p  - 1, \\
   & p(T,x)=0. \nonumber
\end{align}
Here,
$G^{u}$ is a bounded operator from $\calh$ to $\calv'$,
defined by
\begin{align*}
     _{\calv}(\vf, G^{u} \psi)_{\calv'}
     :=\int_{H} \l<Bu(x), Q^\frac 12 D\vf(x)\r> \psi(x) \nu(dx),\  \ \vf \in  \calv,\ \psi\in \calh.
\end{align*}
Note that
$ \|G^{u}\psi\|_{ \calv'} \leq  \|u\|_{L^\9(H;H,\nu)} \|B\| \|\psi\|_{ \calh}$,
where $\|B\|$ denotes the  operator norm.

The backward equation \eqref{equa-p} is understood in
the variational sense,
and its global well-posedness
is part of Theorem \ref{Thm-KE-Back} of Subsection \ref{Subsec-OP-KE-Symm} below.

Now,
we are ready to state the first-order necessary condition for the optimal feedback controllers
in the symmetric case (see Subsection \ref{Subsec-OP-KE-Symm}).

\begin{theorem} \label{Thm-Char-Optimal-Symm}
({\it Necessary condition of optimality: the symmetric case})

Assume that the conditions of Theorem \ref{Thm-Contr-Symm} hold
and let $u_*$ be an optimal controller for Problem $(P^*)$.
Then, we have
\begin{align} \label{u*-character}
  \int_H \l<B(u-u_*), \int_0^T Q^\frac 12 D \vf_* p_* dt \r> d\nu \geq 0,\ \ \forall u\in \calu_{ad},
\end{align}
where
$\vf_*$ and $p_*$
are  the solutions to
\eqref{equa-Nu} and \eqref{equa-p}, respectively,
with $u_*$ replacing $u$.
\end{theorem}

Below
we consider the optimal control problem of the original stochastic differential equation \eqref{equa-X}.

As mentioned in the Introduction,
the concept of martingale problem is robust under bounded perturbations
which, in particular,
fits the optimal control problems considered here.
Moreover,
the martingale problem is well posed in a quite general setting
(e.g., the nonlinearity $F$ may be not continuous
or the operator $Q^\frac 12$
is not necessary Hilbert-Schmidt)
in which case
probabilistic strong solutions may not exist.

Following \cite{BBDR06,DR02},
the martingale problem for \eqref{equa-X} is defined
in Definition \ref{Def-equacontr},
where we use the notion of ``$\nu$-martingale problem''
to express its dependence on the probability measure $\nu$ on $H$.
Later, however,
we shall fix $\nu$ as in Hypothesis $(H1)$ for the remaining part of the paper,
and for simplicity we shall drop the pre-fix $\nu$ again.

\begin{definition} \label{Def-equacontr}
Let $\nu$ be a Borel probability measure on $H$.
A solution to the $\nu$-martingale problem for $(N_0^u, \calf C_b^2(H))$
is a conservative Markov process
$\mathbb{M}^u=(\Omega, \mathscr{F}$, $(\mathscr{F}_t)_{t\geq 0}$, $(X^u(t))_{t\geq 0}, (\bbp_x)_{x\in H_0})$
on $H_0: = {\rm supp} (\nu)$
with  continuous sample paths $t\mapsto\l<X^u(t), e_i\r>$, $i\geq 1$, such that
$X^u(0) =x$, $\bbp_x$-a.s.,
and the following properties hold:
\begin{enumerate}
  \item[(i)] There exist $M, \ve \in (0,\9)$  such that
  \begin{align} \label{def-pt}
      \int_{H_0} (P^u_t f)^2 d\nu \leq M \int_{H_0} f^2 d\nu,\ \
      \forall f\in C_b(H),\ t\in(0,\ve),
  \end{align}
  where $P^u_t$, $t\geq 0$,
  is the transition semigroup of $\bbm^u$.
  \item[(ii)]
  For $\nu$-a.e. $x \in H$,   $\bbp_x$-a.s.,
  \begin{align} \label{A-bdd-Def}
       \int_0^t |\<A(X^u(s)), e_i\>|ds <\9, \ \ for\ every\ t>0,
  \end{align}
  and for all test functions $\vf \in \calf C_b^2(H)$
  \begin{align} \label{vf-X-mart}
       \vf(X^u(t)) - \int_0^t (N^u_0 \vf (X^u(s)) ds,\ \ t\geq 0,
  \end{align}
  is an $(\mathscr{F}_t)$-martingale under $\bbp_x$.
\end{enumerate}
For simplicity,
we also say that $X^u$ solves the
$\nu$-martingale problem for \eqref{equa-X}.

Uniqueness holds if any two
Markov processes which are solutions to the $\nu$-martingale problem for \eqref{equa-X}
have the same finite dimensional distributions
for $\nu$-a.e. starting points $x\in H$.

The $\nu$-martingale problem for \eqref{equa-X} is said to be well posed
if  existence and uniqueness of solutions hold.
\end{definition}

\begin{remark} \label{Rem-A-bdd}
$(i)$ We note that \eqref{A-bdd-Def} holds under Hypothesis $(H1)$,
by the integrability properties of $|x|_H$ and $F$ in Hypothesis $(H1)$ $(i)$
and \eqref{def-pt}.

$(ii)$ The uniqueness of solutions to the martingale problem
can be derived from the existence of martingale solutions
and the m-dissipativity of Kolmogorov operators in Hypothesis $(H1)$.
See the arguments in the proof of the uniqueness part of Theorem \ref{Thm-GWP-Mart}.
\end{remark}

In order to consider the optimal control problem of
the stochastic equation \eqref{equa-X},
we assume that
\begin{enumerate}
  \item[(H1)']  The $\nu$-martingale problem for \eqref{equa-X} is well posed   in the case  $u\equiv 0$.
\end{enumerate}

\begin{remark} \label{Rem-WP-Mart}
$(i)$
Hypothesis $(H1)'$ can be obtained from $(H1)$
if the associated generalized Dirichlet form
is quasi-regular and has the local property.
Actually,
under $(H1)$,
it is known that
(see \cite[p.6]{S99}) the closure $(N_2, D(N_2))$ of $(N_0, \calf C_b^2(H))$
induces a generalized Dirichlet form on $L^2(H,\nu)$
as follows
\begin{align*}
   \cale(\vf, \psi):= \left\{
                        \begin{array}{ll}
                           - (N_2\vf, \psi), & \hbox{$\vf \in D(N_2)$, $\psi \in L^2(H,\nu)$;} \\
                           - (N_2^* \psi, \vf), & \hbox{$\vf \in L^2(H,\nu)$, $\psi \in D(N_2^*)$,}
                        \end{array}
                      \right.
\end{align*}
where $(\ ,\ )$ denotes the inner product in $L^2(H,\nu)$,
and $N_2^*$ is the dual operator of $N_2$.
We see that the condition $D3$ in \cite[p.78]{S99} is satisfied
with $\mathcal{Y}=\calf C_b^2(H)$ and $\calf = D(N_2)$.
Hence,
provided $\cale$ is quasi-regular,
\cite[Theorem IV. 2.2]{S99} yields that
there exists a sufficiently regular Markov process $\bbm$
(namely, a $\nu$-tight special standard process)
with transition semigroup $P_t$, $t>0$,
given by $P_t^\nu$, $t>0$,
hence satisfying $(i)$ of Definition \ref{Def-equacontr}
with $M=1$ for all $t\in (0,\9)$.
In addition,
by \cite[Theorem 3.3]{T03},
the sample paths of $\bbm$ are continuous $\bbp$-a.s.
for $\nu$-a.e. $x\in H_0$
if $\mathcal{E}$ is local.
Moreover,
Remark \ref{Rem-A-bdd} $(i)$ and \cite[Proposition 1.4]{BBR06} yield that
$\bbm$ satisfies the property $(ii)$ of Definition \ref{Def-equacontr}.
Thus, $\bbm$ solves the martingale problem for \eqref{equa-X} for $u\equiv 0$.
The uniqueness can be proved as in Section \ref{Sec-OP-SDE} below.

As specific examples,
we show in Theorems \ref{Thm-H1iv} and \ref{Thm-H1'-H1-Sym} that
Hypothesis $(H1)'$ can be implied by the m-dissipativity of Kolmogorov operators
in certain situations,
based on the theory of (generalized) Dirichlet forms.

$(ii)$
It follows by \eqref{vf-X-mart} and $(H1)$ that for all $\vf\in \calf C_b^2(H)$
and for all $t>0$
\begin{align*}
     H_0 \ni x \mapsto P^u_t \vf(x):= \bbe_{\bbp_x} [\vf(X^u(t))]
\end{align*}
is a $\nu$-version of $e^{tN_2^u}\vf \in L^2(H, \nu)$.
Below we shall briefly describe this by saying that
``$P^u_t$ is given by $e^{tN_2^u}$''.
\end{remark}

Similarly to Theorem \ref{Thm-Max-Nu},
the well-posedness of the martingale problem for controlled equations
can be inherited from that for the uncontrolled equation.

\begin{theorem} \label{Thm-GWP-Mart}
Assume that Hypotheses $(H1)$ and $(H1)'$ hold.
Then, for each $u\in \calu_{ad}$,
the martingale problem for \eqref{equa-X} is well posed.
\end{theorem}

The main result for optimal control problems of
the stochastic equation \eqref{equa-X} is formulated below.

\begin{theorem} \label{Thm-Contr-SDE}
({\it Optimal control for stochastic equations on Hilbert spaces})

$(i)$ ({\it General case}) Assume $(H1)$ and $(H1)'$.
Assume in addition that either $(H2)$ or $(H3)$ holds.
Then, for any objective function $g\in D(N_2)$,
there exists an optimal control to Problem $(P)$.

$(ii)$ ({\it Symmetric case})
Consider the symmetric case.
Assume $(H1)$ and that
$(Q^\frac 12 D$, $\mathcal{F}C_b^1(H))$ is closable from $L^2(H,\nu)$ to $L^2(H;H,\nu)$.
Assume additionally $(H2)$ or $(H3)$.
Then, for any objective function $g\in L^2(H,\nu)$,
there exists an optimal control to Problem $(P)$
and the first-order necessary condition \eqref{u*-character} holds.

In both cases, when $(H3)$ holds, one can take $B=Id$.
\end{theorem}

The remainder of this paper is organized as follows.
Section \ref{Sec-OP-KE}
is mainly devoted to the optimal control problem of Kolmogorov equations.
We first prove Theorem \ref{Thm-Contr} in the general case in Subsection \ref{Subsec-OP-KE-General},
while  Theorems \ref{Thm-Contr-Symm} and \ref{Thm-Char-Optimal-Symm}
are proved later in Subsection \ref{Subsec-OP-KE-Symm}.
Then, in Section \ref{Sec-OP-SDE}
we study the optimal control problem of
stochastic equations on Hilbert spaces.
Finally, Section \ref{Sec-Application} mainly
contains the applications to
various stochastic partial differential equations,
including  stochastic  equations in $H$ with singular drifts,
stochastic reaction-diffusion equations
and stochastic porous media equations.

\section{Optimal control of Kolmogorov equations} \label{Sec-OP-KE}

\subsection{General case} \label{Subsec-OP-KE-General}

This subsection is mainly devoted to the proof of Theorem \ref{Thm-Contr}.
To begin with,
we first prove Theorem \ref{Thm-Max-Nu}
for the realted Kolmogorov operators.\\

{\bf Proof of Theorem  \ref{Thm-Max-Nu}.}
For simplicity, we set $\calh:= L^2(H,\nu)$.
Let us first prove the identity \eqref{Integ-part}.
Actually, by
straightforward computations,
\begin{align*}
   N_0(\vf^2) = 2 \vf N_0 \vf + |Q^\frac 12 D \vf|_H^2, \ \ \forall \vf \in \calf C_b^2(H),
\end{align*}
which along with  $(H1)$ $(ii)$ implies that
\begin{align}   \label{integ-part-N0-2}
   \int_H  \vf N_0 \vf d\nu
   = -\frac 12 \int_H |Q^\frac 12 D \vf|_H^2 d\nu ,  \ \ \forall \vf \in \calf C_b^2(H).
\end{align}
To extend \eqref{integ-part-N0-2} to all $\vf \in D(N_2)$
we observe that the map
\begin{align*}
   Q^\frac 12 D: \calf C_b^2(H) \to L^2(H;H,\nu)
\end{align*}
is linear and by \eqref{integ-part-N0-2}
continuous with respect to the graph norm of $N_2$ on $\calf C_b^2(H)$.
Since $\calf C_b^2(H)$ is dense in $D(N_2)$
with respect to the graph norm of $N_2$,
this map extends uniquely by continuity to $D(N_2)$ and
then \eqref{Integ-part} follows by continuity.

We also note that
\begin{align} \label{perb-def}
    D(N_2) \ni \vf \mapsto \langle Bu, \ol{Q^\frac 12 D} \vf\rangle \in \calh
\end{align}
is a well-defined bounded linear operator.
Actually, for any $\vf \in D(N_2)$,
we take $\{\vf_n\}\subseteq \calf C_b^2(H)$ such that
$\vf_n \to \vf$ in the $N_2$-graph norm in $\calh$.
Then,
since $\sup_{x\in H} |u(x)|_H \leq \rho$,
by \eqref{integ-part-N0-2},
\begin{align} \label{esti-perb*}
    \int |\l< Bu, \ol{Q^\frac 12 D} (\vf_n -\vf_m) \r>|^2 d\nu
    \leq& \rho^2 \|B\|^2 \int |\ol{Q^\frac 12 D}(\vf_n-\vf_m)|^2 d\nu \\
    =& -2 \rho^2 \|B\|^2 \int  (\vf_n -\vf_m)  N_0 (\vf_n -\vf_m)d\nu
    \to 0,   \nonumber
\end{align}
as $n,m\to \9$.
This implies that $\{\l<Bu, \ol{Q^\frac 12 D} \vf_n\r>\}$ is a Cauchy sequence in $\calh$
and so yields the claim above.

In order to prove that $\calf C_b^2(H)$ is a core of $(N_2^u, D(N_2))$,
it suffices to prove that
the graph norms of $N_2$ and $N_2^u$ are equivalent, i.e.,
there exists $C>0$ such that for any $\vf \in \mathcal{F} C_b^2(H)$,
\begin{align} \label{equi-N2-N2u-norm}
     C^{-1} (\|N_2 \vf\|_\calh^2 + \|\vf\|_\calh^2)
     \leq (\|N^u_2 \vf\|_\calh^2 + \|\vf\|_\calh^2)
     \leq C (\|N_2 \vf\|_\calh^2 + \|\vf\|_\calh^2).
\end{align}
For this purpose,
we note that for any $\lbb >0$,
similarly to \eqref{esti-perb*},
\begin{align} \label{esti-perb.1}
   \| \l<B u,  \ol{Q^\frac 12 D} \vf \r>\|^2_\calh
  \leq  -2  \rho^2 \|B\|^2 \int_H    \vf  N_2 \vf d\nu
   \leq  2\rho^2 \|B\|^2 (\frac 1 \lbb \|N_2\vf\|_\calh^2
         +   \lbb \|\vf\|^2_\calh),
\end{align}
which immediately yields the second inequality of \eqref{equi-N2-N2u-norm}.
Moreover,
taking $\lbb$ large enough such that $2\rho^2 \|B\|^2/\lbb \leq 1/4$
we also obtain the first inequality of \eqref{equi-N2-N2u-norm}.

The fact that $(N_2^u, D(N_2))$ generates a $C_0$-semigroup $e^{tN_2^u}$
follows from the essential m-dissipativity of $(N_0^u, \mathcal{F}C_b^2(H))$
on $\calh$.
To this end,
we note that
for $\lbb$ large enough,
for any $f\in \calh $,
the equation
\begin{align*}
   \lbb \vf - N_2 \vf - \l<B u, \ol{Q^\frac 12 D} \vf \r> =f
\end{align*}
has the unique solution
\begin{align*}
   \vf = R_\lbb ((I- T_\lbb)^{-1} f),
\end{align*}
where $R_\lbb$ is the resolvent of $N_2$, i.e.,
$R_\lbb = (\lbb - N_2)^{-1}$,
and the operator
$T_\lbb: L^2(H, \nu) \to L^2(H, \nu)$
is defined by $T_\lbb \psi = \l<B u, \ol{Q^\frac 12 D} R_\lbb \psi\r>$.
(Note that, $\|T_\lbb \psi\|_{\calh} \leq \frac 12 \|\psi\|_\calh$
when $\lbb$ is large enough,
hence $(I- T_\lbb)^{-1}$ is well-defined
and $(I-T_\lbb)^{-1} \in L(\calh)$.)
By the essential m-dissipativity of $(N_0, \calf C_b^2(H))$,
there exists a sequence $(\vf_n) \subseteq \calf C_b^2(H)$
such that
$(\lbb - N_2) \vf_n \to (I-T_\lbb)^{-1} f$, as $n\to \9$.
This yields
\begin{align*}
   (\lbb - N_2) \vf_n - \l<B u(x), \ol{Q^\frac 12 D}\vf_n\r>
   =&  (\lbb - N_2) \vf_n - T_\lbb (\lbb - N_2)\vf_n  \\
   \to & (I-T_\lbb)^{-1} f - T_\lbb (I-T_\lbb)^{-1} f
   = f,
\end{align*}
as $n\to \9$,
which implies that
the image of $\lbb - N_0^u$ is dense in $\calh$.
Thus,
$(N_0^u, \mathcal{F}C_b^2(H))$ is essentially m-dissipative on $\calh$
and so generates a semigroup $e^{tN_2^u}$ on $\calh$, $t\geq 0$.

Regarding
\eqref{bdd-eN2u},
by \eqref{Integ-part} and Cauchy's inequality,
for $\vf := e^{t N^u_2}g$,
\begin{align*}
   \frac 12 \p_t \|\vf\|_\calh^2
   =& (N_2 \vf, \vf) + (\l<Bu, \ol{Q^\frac 12 D} \vf\r>, \vf) \\
   \leq& -\frac 12 \int_H |\ol{Q^\frac 12 D} \vf|_H^2 d \nu
         + \rho \|B\| \int_H |\ol{Q^\frac 12 D} \vf|_H |\vf| d \nu  \\
   \leq& -\frac 14  \int_H |\ol{Q^\frac 12 D} \vf|_H^2 d \nu
          + 4 \rho^2 \|B\|^2 \|\vf\|_\calh^2,
\end{align*}
which along with Gronwall's inequality implies \eqref{bdd-eN2u}.

Therefore, the proof of Theorem \ref{Thm-Max-Nu} is complete.
\hfill $\square$\\

{\bf Proof of Theorem \ref{Thm-Contr}.}
We set $\calh := L^2(H, \nu)$ and
\begin{align*}
   \Phi(u): = \int_0^T \int_H \psi^u(t,x) \nu(dx)dt,
\end{align*}
where $\psi^u$ is the solution to \eqref{equa-Nu}
corresponding to $u\in \calu_{ad}$.

Let $I_* := \inf\{\Phi(u): u\in \calu_{ad}\}$
and take a sequence $\{u_n\}  \subseteq\calu_{ad}$, such that
$I_*\leq \Phi(u_n) \leq I_* +\frac 1n$, $n\geq 1$.

Let  $\vf_n  := e^{- (4 \rho^2 \|B\|^2+1) t} \psi^{u_n}$,
$n\geq 1$.
Then, we have
\begin{align} \label{Phi-un}
   \Phi(u_n) = \int_0^T \int_H e^{(4 \rho^2 \|B\|^2+1)t} \vf_n d \nu  dt,
\end{align}
and by \eqref{equa-Nu} $\vf_n$ solves
\begin{align} \label{equa-calN}
     & \frac{d\vf_n}{dt} = \caln^{u_n} \vf_n ,\ \ t\in (0,T),\\
     & \vf_n(0) = g, \nonumber
\end{align}
where  the operator $\caln^{u_n}: D(N_2) \to \calh$ is defined by
\begin{align} \label{calnun}
  \caln^{u_n} \psi:= N_2 \psi + \l<B u_n, \ol{Q^\frac 12 D} \psi\r> - (4 \rho^2 \|B\|^2+1) \psi, \ \ \psi \in D(N_2).
\end{align}

By Theorem \ref{Thm-Max-Nu},
we see that also
$(\caln^{u_n}, D(N_2))$ generates a $C_0$-semigroup $e^{t\caln^{u_n}}$ on $\calh$,
namely $e^{t\caln^{u_n}} = e^{-(4 \rho^2 \|B\|^2+1) t} e^{tN_2^{u_n}}$,
where $ e^{tN_2^{u_n}}$ is given by Theorem \ref{Thm-Max-Nu},
and
\begin{align} \label{esti-vfj}
    \sup\limits_{n\geq 1}
      \|\vf_n\|^2_{C([0,T];\calh)}
    + \sup\limits_{n\geq 1}
      \int_0^T \int |\ol{ Q^\frac 12 D} \vf_n(s)|_H^2 d \nu   ds
    \leq C \|g\|^2_{\calh}.
\end{align}
Similarly, by \eqref{equi-N2-N2u-norm} and \eqref{esti-vfj},
\begin{align} \label{bdd-Nunvfn}
  \sup\limits_{n\geq 1} \|\caln^{u_n}\vf_n\|_{C[0,T; \calh]}
   = \sup\limits_{n\geq 1}  \|e^{t\caln^{u_n}} \caln^{u_n} g\|_{C[0,T; \calh]}
   \leq  C(T,\rho) (\|N_2g\|_\calh + \|g\|_\calh).
\end{align}

Hence, along a subsequence,
again denoted $\{n\}$, we have
\begin{align}
    u_n \to u_*,\ \     &weakly-star\ in\ L^\9(H; H, \nu), \label{uj-weakH-rd} \\
    \vf_n \to \vf_*,\ \ &weakly\ in\ L^2(0,T; \calh), \label{vfj-weakH-rd} \\
    \caln^{u_n} \vf_n\to \eta,\ \ &weakly\ in\ L^2(0,T; \calh).    \label{Nuj-weakV-rd}
\end{align}

Note that,
since $t \mapsto e^{(4 \rho^2 \|B\|^2+1)t} \in L^2(0,T; \calh)$,
we apply \eqref{vfj-weakH-rd} to pass to the limit in \eqref{Phi-un}
to obtain
\begin{align} \label{Phi-I*}
   I_* = \lim\limits_{n\to \9}  \Phi(u_n)
   = \int_0^T \int_H e^{4\rho^2  \|B\|^2 t} \vf_* d \nu dt.
\end{align}
Moreover, by \eqref{equa-calN}, \eqref{vfj-weakH-rd} and \eqref{Nuj-weakV-rd},
$\vf_*$ satisfies the equation
\begin{align} \label{equa-rd-eta}
    \vf_*(t) = g+   \int_0^t \eta(s) ds,\ \ for\ dt-a.e.\ t\in [0,T].
\end{align}

Now it remains to prove that
\begin{align} \label{eta-uvf-rd}
     \eta (t) = \caln^{u_*} \vf_*(t),\ for\ a.e.\ t\in (0,T).
\end{align}

Note that, by \eqref{Integ-part},
\begin{align} \label{calnu-mono}
    (\caln^{u_n} v, v)_{\calh} \leq 0,\ \ \forall v\in D(N_2),
\end{align}
For simplicity,
we set $\calh_t:= L^2(0,t;\calh)$
with the inner produce $(\cdot, \cdot)_{\calh_t}$ below.
Then, \eqref{calnu-mono}  yields that, for any positive function $h\in L^\9(0,T)$
and any $\psi \in \mathcal{M}$,
where $ \mathcal{M}$ denotes the space of all linear combinations of functions of the form
$fv$, where $f\in L^\9(0,T)$ and $v\in \calf C_b^2(H)$,
\begin{align} \label{0-K1234}
    0 \geq& \int_0^T h(t) (\caln^{u_n} (\vf_n - \psi), \vf_n - \psi)_{\calh_t} dt  \nonumber \\
      =& \int_0^T h(t)  (\caln^{u_n} \vf_n, \vf_n)_{\calh_t} dt
         - \int_0^T h(t) (\caln^{u_n} \vf_n, \psi)_{\calh_t} dt   \nonumber \\
      & - \int_0^T h(t) (\caln^{u_n} \psi, \vf_n)_{\calh_t} dt
         + \int_0^T h(t) (\caln^{u_n} \psi, \psi)_{\calh_t} dt  \nonumber \\
      =:& K_{1,n} - K_{2,n} -  K_{3,n} + K_{4,n}.
\end{align}

Below we treat $K_{i,n}$ separately, $1\leq i\leq 4$.

For $K_{1,n}$, by equation \eqref{equa-calN},
\begin{align*}
   \frac 12 \|\vf_n(t)\|^2_\calh = (\caln^{u_n} \vf_n, \vf_n )_{\calh_t} + \frac 12 \|g\|_\calh^2, \ \ 0<t<T.
\end{align*}
Then, multiplying both sides by $h(t)$ and integrating over $[0,T]$ we get
\begin{align} \label{esti-K1j.1}
    \frac 12 \int_0^T h(t) \|\vf_n(t)\|^2_\calh dt
    = \int_0^T h(t) (\caln^{u_n} \vf_n, \vf_n )_{\calh_t} dt
      + \frac 12 \|g\|_\calh^2 \int_0^T h(t) dt.
\end{align}
Similarly,
we infer from the equation \eqref{equa-rd-eta} that
\begin{align} \label{esti-K1j.2}
    \frac 12 \int_0^T h(t) \|\vf_*(t)\|^2_\calh dt
    = \int_0^T h(t)  (\eta, \vf_*)_{\calh_t} dt
      + \frac 12 \|g\|_\calh^2 \int_0^T h(t) dt.
\end{align}
Moreover, by \eqref{vfj-weakH-rd},
\begin{align} \label{weak-vfnh}
     \vf_n h^\frac 12  \to \vf_* h^\frac 12,\ \ weakly\ in\ L^2(0,T; \calh),
\end{align}
which implies that
\begin{align} \label{esti-K1j.3}
     \liminf\limits_{n\to \9}
     \int_0^T h(t) \|\vf_n(t)\|_\calh^2 dt
     \geq \int_0^T h(t) \|\vf_*(t)\|_\calh^2 dt.
\end{align}
Thus, we obtain from \eqref{esti-K1j.1}-\eqref{esti-K1j.3} that
\begin{align} \label{lim-K1j-rd}
     \liminf\limits_{n\to \9}  K_{1,n}
     \geq \int_0^T h(t)  (\eta, \vf_*)_{\calh_t} dt.
\end{align}

Moreover, in order to pass to the limit in $K_{2,n}$,
we note that, by \eqref{Nuj-weakV-rd},
\begin{align*}
     (\caln^{u_n} \vf_n, \psi)_{\calh_t}
     \to (\eta, \psi)_{\calh_t}, \ \ \forall t\in (0,T].
\end{align*}
Taking into account \eqref{bdd-Nunvfn},
we have for any $t\in (0,T]$,
\begin{align*}
     | (\caln^{u_n} \vf_n, \psi)_{\calh_t} |
     \leq & \|\caln^{u_n} \vf_n\|_{L^2(0,T; \calh)} \|\psi\|_{L^2(0,T; \calh)} \\
     \leq & C (T,\rho) (\|N_2 g\|_\calh + \|g\|_\calh) \|\psi\|_{L^2(0,T; \calh)} <\9.
\end{align*}
Then, the dominated convergence theorem yields
\begin{align} \label{lim-K2j-rd}
     \lim\limits_{n\to \9}  K_{2,n}
     = \int_0^T h(t) (\eta, \psi)_{\calh_t} dt.
\end{align}

Now, we treat the  term $K_{3,n}$.
We expand
\begin{align} \label{lim-K3j.0}
   (\caln^{u_n} \psi, \vf_n)_{\calh_t}
   =  (N_2 \psi, \vf_n)_{\calh_t}
      +  ( \l< B u_n, \ol{Q^\frac 12 D}\psi\r>, \vf_n)_{\calh_t}
      - 4\rho^2 \|B\|^2 (\psi, \vf_n)_{\calh_t}.
\end{align}
By \eqref{vfj-weakH-rd},
\begin{align} \label{lim-K3j.13}
     \lim\limits_{n\to \9} (N_2 \psi, \vf_n)_{\calh_t}
     =  (N_2 \psi, \vf_*)_{\calh_t}, \ \
      \lim\limits_{n\to \9}  ( \psi, \vf_n)_{\calh_t}
     =  ( \psi, \vf_*)_{\calh_t}.
\end{align}
Concerning the second term on the right-hand side of \eqref{lim-K3j.0},
we claim that
\begin{align} \label{lim-K3j.2}
     ( \l< B u_n, \ol{Q^\frac 12 D} \psi\r>, \vf_n)_{\calh_t}
    \to  ( \l< B u_*, \ol{Q^\frac 12 D} \psi\r>, \vf_*)_{\calh_t}.
\end{align}

In order to prove \eqref{lim-K3j.2}, we note that
\begin{align} \label{lim-K3j.2*}
   &   ( \l< B u_n, \ol{Q^\frac 12 D} \psi\r>, \vf_n)_{\calh_t}
     -  ( \l< B u_*, \ol{Q^\frac 12 D} \psi\r>, \vf_*)_{\calh_t}   \\
   =&   ( \l< B u_n- B u_*, \ol{Q^\frac 12 D} \psi\r>, \vf_*)_{\calh_t}
          +  \int_0^t \int_H \l< B u_n, (\vf_n(s) - \vf_*(s)) \ol{Q^\frac 12 D} \psi(s)\r> d\nu ds. \nonumber
\end{align}
For the first term
on the right-hand side of \eqref{lim-K3j.2*},
Fubini's theorem yields
\begin{align*}
      (  \l< B u_n- B u_*, \ol{Q^\frac 12 D} \psi\r>, \vf_*)_{\calh_t}
   =& \int_0^t \int_H \l< Bu_n -  Bu_*, \vf_*(s) \ol{Q^\frac 12 D} \psi(s)\r> d\nu ds \\
   =& \int_H \l<Bu_n -  B u_*, \int_0^t \vf_*(s) \ol{Q^\frac 12 D} \psi (s) ds\r> d\nu.
\end{align*}
But the latter term converges to zero,
since $ B \in L(L^\9(H;H,\nu))$
and thus $ Bu_n \to B u_*$ as $n\to \9$
weakly-star in $L^\9(H;H,\nu)$
and since for all $t\in [0,T]$,
\begin{align*}
   \int_0^t \vf_*(s) \ol{Q^\frac 12 D} \psi(s) ds \in L^1 (H; H,\nu).
\end{align*}

Now, let us treat the second term on the right-hand side of \eqref{lim-K3j.2*}.
If $(H2)$ holds, then
it follows that $ Bu_n \to B u_*$ in $L^2(H;H,\nu)$ as $n\to \9$,
hence the second term in \eqref{lim-K3j.2*} obviously converges to zero as $n\to \9$.

Now assume that $(H3)$ holds and set $\calv : = W^{1,2}(H,\nu)$.
We  claim that
\begin{align} \label{vfj-stroH-rd}
    \vf_n \to \vf_*,\ \ strongly\ in\ L^2(0,T; \calh).
\end{align}
Actually,  from equation \eqref{equa-Nu}
it follows that $\frac{d\vf_n}{dt}$ is bounded in $L^2(0,T; \calh)$,
and by \eqref{esti-vfj},
$\{\vf_n\}$ is also bounded in $L^2(0,T; \calv)$.
Since by Hypothesis $(H3)$,
$\calv$ is compactly embedded in $\calh$,
by virtue of the Aubin-Lions compactness theorem
(see e.g. \cite[Theorem 1.2]{B10})
we obtain \eqref{vfj-stroH-rd}, as claimed.

Hence, the second term on the right-hand side of \eqref{lim-K3j.2*} also converges if $(H3)$ holds.

Thus,
we obtain \eqref{lim-K3j.2}, as claimed.
(Note that, Hypothesis $(H2)$ or $(H3)$
are only used here.
When $(H3)$ holds one may take $B=Id$.)

This along with  \eqref{lim-K3j.0} and \eqref{lim-K3j.13} yields
\begin{align} \label{Nujvf-Nuvf}
     \lim\limits_{n\to \9} {}
      (\caln^{u_n} \psi, \vf_n)_{\calh_t}
     =&  (N_2 \psi, \vf_*)_{\calh_t}
       +  (\l<  B u_*, \ol{Q^\frac 12 D} \psi\r>, \vf_*)_{\calh_t}
       -4\rho^2 \|B\|^2  (\psi, \vf_*)_{\calh_t}  \nonumber \\
     =&  (\caln^{u_*} \psi, \vf_*)_{\calh_t}.
\end{align}
Taking into account that
\begin{align*}
     \sup\limits_{t\in [0,T]}
     | (\caln^{u_n} \psi, \vf_n)_{\calh_t} |
     \leq  C  \sup\limits_{j\geq 1}  \|\vf_n\|_{L^2(0,T; \calh)}<\9,
\end{align*}
we may apply the  dominated convergence theorem to obtain
\begin{align} \label{lim-K3j-rd}
     \lim\limits_{n\to \9} K_{3,n}
     = \int_0^T h(t)   (\caln^{u_*} \psi, \vf_*)_{\calh_t} dt.
\end{align}

Finally, arguing as in the proof of  \eqref{lim-K3j-rd} we have
\begin{align} \label{lim-K4j-rd}
     \lim\limits_{n\to \9} K_{4,n}
     = \int_0^T h(t)   (\caln^{u_*} \psi, \psi)_{\calh_t} dt.
\end{align}

Thus, combing together
\eqref{0-K1234}, \eqref{lim-K1j-rd}, \eqref{lim-K2j-rd}, \eqref{lim-K3j-rd} and \eqref{lim-K4j-rd}
we conclude that
for any positive functions $h\in L^\9(0,T)$
and any $\psi \in \mathcal{M}$,
\begin{align} \label{0-eta-caln-rd}
    0\geq \int_0^T h(t)
          (\eta - \caln^{u_*} \psi, \vf_*- \psi)_{\calh_t} dt.
\end{align}

Now, consider the Hilbert space  $L^2(0,T; D(N_2))$,
where $D(N_2)$ is equipped with the graph norm given by $\caln^{u_*}$.
Then, if $\psi \in L^2(0,T; D(N_2))$
such that
\begin{align*}
   \int_0^T \<\psi, fv\>_{D(N_2)} dt = 0,\ \ \forall f\in L^\9(0,T),\ v\in \calf C_b^2(H),
\end{align*}
since $\calf C_b^2(H)$ is dense in $D(N_2)$,
it follows that $\psi =0$,
i.e., $\mathcal{M}$ is dense in $L^2(0,T; D(N_2))$.

Hence, for every $\psi \in L^2(0,T; D(N_2))$,
there exist $\psi_n \in \mathcal{M}$,
such that $\psi_n \to \psi$ in $L^2(0,T; D(N_2))$ as $n\to \9$.
Then, we have
\begin{align*}
  (\eta - \caln^{u_*} \psi, \vf_* - \psi)_{\calh_t}
  = \lim\limits_{n\to \9}  (\eta - \caln^{u_*} \psi_n, \vf_* - \psi_n)_{\calh_t}
\end{align*}
and for all $t\in [0,T]$,
\begin{align*}
     |(\eta - \caln^{u_*} \psi_n, \vf_* - \psi_n)_{\calh_t} |
  \leq& \|\eta\|_{L^2(0,T; \calh)}(\|\vf_*\|_{L^2(0,T; \calh)} + \sup\limits_{n\geq 1} \| \psi_n\|_{L^2(0,T; \calh)})  \\
      & + \sup\limits_{n\geq 1} \|\caln^{u_*} \psi_n\|_{L^2(0,T; \calh)}  (\|\vf_*\|_{L^2(0,T; \calh)} + \sup\limits_{n\geq 1} \| \psi_n\|_{L^2(0,T; \calh)})
\end{align*}

Hence, \eqref{0-eta-caln-rd} extends to all $\psi \in L^2(0,T; D(N_2))$
and, if $\vf_* \in L^2(0,T; D(N_2))$,
we may particularly take $\psi := \vf_*(t) - \ve f(t) v$,
where $\ve \in (0,1)$, $f \in L^\9(0,T)$ and $v\in \calf C_b^2(H)$.
Dividing by $\ve$ and then letting $\ve$ to $0$,
we arrive at
\begin{align} \label{0-eta-calnu}
     0 \geq \int_0^T h(t)
           (\eta - \caln^{u_*} \vf_*, f v)_{\calh_t} dt.
\end{align}

Therefore,
since $h,f \in L^\9(0,T)$ and $v\in \calf C_b^2(H)$
are arbitrary,
we obtain \eqref{eta-uvf-rd},
if $\vf_* \in L^2(0,T; D(N_2))$.

To show that $\vf_* \in L^2(0,T; D(N_2))$,
we first note that
because of \eqref{vfj-weakH-rd} we only have to prove that
\begin{align}   \label{calnun-vfn-bdd}
   \sup\limits_{n\geq 1} \|\caln^{u_*} \vf_n\|_{L^2(0,T; \calh)} <\9.
\end{align}
Indeed, then there exists a subsequence of $\vf_n$, $n\in \mathbb{N}$,
along which $\caln^{u_*} \vf_n$, $n\in \mathbb{N}$,
converges weakly in $L^2(0,T; \calh)$.
Hence,
taking into account \eqref{vfj-weakH-rd}
we can select a subsequence $\vf_{n_l}$, $l\in \mathbb{N}$,
such that the following two Cesaro-mean convergence strongly in $L^2(0,T; \calh)$
(by the Banach-Saks Theorem),
\begin{align} \label{vf-strong-conv}
    \frac 1 N \sum\limits_{l=1}^N \vf_{n_l} \to \vf_*,\ \ in\ L^2(0,T; \calh),\ as\ N\to\9,
\end{align}
and
\begin{align} \label{Nuvf-strong-conv}
   \caln^{u_*} (\frac 1 N \sum_{l=1}^N \vf_{n_l}),\ \ N\in \mathbb{N},\ converges\  in\ L^2(0,T; \calh),
\end{align}
so $\vf_* \in L^2(0,T; D(N_2))$ by completeness.

But we have by \eqref{calnun},  \eqref{esti-vfj}, \eqref{Nuj-weakV-rd}
and since $\sup_{n\geq 1}\|u_n\|_{L^\9(H; H, \nu)} \leq \rho$, that
\begin{align*}
    \sup\limits_{n\geq 1} \|N_2 \vf_n\|_{L^2(0,T,\calh)} <\9.
\end{align*}
Hence, it follows by \eqref{esti-vfj} and  \eqref{calnun} again
that \eqref{calnun-vfn-bdd} holds.

Therefore,
the proof of Theorem \ref{Thm-Contr} is complete.
\hfill $\square$

\subsection{Symmetric case}  \label{Subsec-OP-KE-Symm}

Throughout this subsection, we assume
the self-adjointness of the Kolmogorov operator $N_2$
and the closability of the operator $Q^\frac 12 D$.
Precisely,
we assume that  $N^*_2 = N_2$ on $\calh:= L^2(H,\nu)$,
and the operator $Q^\frac 12 D$ with domain $\mathcal{F} C_b^1(H)$
is closable from $L^2(H,\nu)$ to $L^2(H;H,\nu)$.

A nice feature in the symmetric case is that
the Kolmogorov operators are variational from $\calv$ to $\calv'$.
This enables us to analyze the Kolmogorov equation \eqref{equa-Nu}
and the backward equation \eqref{equa-p}
in the variational setting.

We first extend the integration by parts formula \eqref{Integ-part} in the symmetric case.
\begin{lemma} \label{Lem-integ-symm}
Consider the symmetric case.
Assume $(H1)$
and that
$(Q^\frac 12 D, \calf C_b^1(H))$ is closable from $L^2(H, \nu)$ to $L^2(H; H, \nu)$.
Then,
\begin{align} \label{integ-part-symm}
   \int_H \psi \ol{N}_2 \vf  d \nu
   = -\frac 12 \int_H \l<Q^\frac 12 D\vf, Q^\frac 12 D \psi\r> d \nu,\ \ \vf,\psi\in \calv:=W^{1,2}(H,\nu),
\end{align}
where $\ol{N}_2$ is the continuous extension of the operator
\begin{align} \label{def-olN2}
     \calf C_b^2(H) \ni \vf \mapsto N_2 \vf \in \calv'
\end{align}
with respect to the $\calv$-norm on $\calf C_b^2(H)$.
\end{lemma}

{\bf Proof.}
Since $N_2^* = N_2$,
using \eqref{integ-part-N0-2} and polarization
we see that \eqref{integ-part-symm} is valid for
all $\vf, \psi \in \calf C_b^2(H)$.
Note that,
by \eqref{integ-part-symm},
the map \eqref{def-olN2} is linear and continuous with respect to the $\calv$-norm
on $\calf C_b^2(H)$.
Hence,
since $\calf C_b^2(H)$ is dense in $\calv$,
the map \eqref{def-olN2} can be extended uniquely by continuity
to $\calv$ and \eqref{integ-part-symm} follows by continuity.
\hfill $\square$ \\

Below we still denote $\ol{N}_2$ by $N_2$ in the symmetric case for simplicity.

In order to obtain the global well-posedness for the Kolmogorov equation \eqref{equa-Nu},
we use $\wt\vf:= e^{-(4\rho^2  \|B\|^2 +1) t} \vf$
to rewrite \eqref{equa-Nu} as follows
\begin{align} \label{equa-wtvf}
     &\frac{d \wt \vf}{dt} = \caln^u \wt\vf, \\
     & \wt\vf(0) = g,  \nonumber
\end{align}
where $\caln^u$ is as in \eqref{calnun},
i.e.,  $\caln^u = N_2^u - (4\rho^2 \|B\|^2 +1)$.

Similarly, for $\wt p(t):= e^{-(4\rho^2 \|B\|^2 +1) t} p (T-t)$,
we have from \eqref{equa-p} that
\begin{align} \label{equa-wtp}
   & \frac{d \wt{p}}{dt} = \wt{\caln}^{u} \wt{p} + e^{- (4\rho^2 \|B\|^2 +1)  t}, \\
   & \wt{p}(0)=0, \nonumber
\end{align}
where $\wt{\caln}^u := N_2 + G^u - (4\rho^2 \|B\|^2 +1)$
with $G^u $ as in \eqref{equa-p}.

The properties of operators $\caln^u$ and $\wt{\caln}^u$
are collected in the result below.

\begin{proposition} \label{Prop-caln-cala}
Under the conditions of Lemma  \ref{Lem-integ-symm},
we have
\begin{align} \label{calna-v-v'}
   \sup\limits_{u\in \calu_{ad}}
   (\|\caln^u \vf\|_{\calv'} + \|\wt{\caln}^u \vf\|_{\calv'}) \leq C(T,\rho) \|\vf\|_{\calv}, \ \ \forall \vf \in \calv,
\end{align}
for some $C(T,\rho)>0$,
and for any $\vf \in \calv$,
\begin{align} \label{calna-dissp}
    _{\calv}(\vf, \caln^u \vf)_{\calv'} + _{\calv}(\vf, \wt{\caln}^u \vf)_{\calv'}
    \leq -\frac 12 \|\vf\|^2_\calv.
\end{align}
\end{proposition}

{\bf Proof.}
Let us first consider the operator $\caln^u$.
By  \eqref{integ-part-symm},
for any $\vf, \psi \in \calv$,
\begin{align} \label{psi-calnuvf}
    _{\calv} \<\psi, \caln^u \vf \>_{\calv'}
   =& -\frac 12 \int_H \l< Q^\frac 12 D \vf, Q^\frac 12 D \psi\r> d \nu
     + \int_H \l<B u, Q^\frac 12 D \vf\r> \psi d \nu \nonumber  \\
    & - 4 \rho^2  \|B\|^2 \int \vf \psi d \nu,
\end{align}
which along with H\"older's inequality implies immediately that
for some $C>0$
\begin{align}
   \| \caln^u \vf \|_{\calv'} \leq C \|\vf\|_{\calv}.
\end{align}
Moreover, by Cauchy's inequality and $ab\leq a^2 + b^2$,
we have
\begin{align*}
   |\int_H \l<B u, Q^\frac 12 D \vf\r> \vf d \nu|
   \leq  \rho \|B\| \int_H |Q^\frac 12 D \vf|_H |\vf| d\nu
   \leq  \frac 14 \|\vf\|_\calv^2 + 4\rho^2 \|B\|^2 \|\vf\|_\calh^2.
\end{align*}
Plugging this into \eqref{psi-calnuvf} with $\vf$ replacing $\psi$ we obtain
\begin{align}
   {}_{\calv} (\vf, \caln^u \vf)_{\calv'}
   \leq - \frac 14 \|\vf\|_\calv^2.
\end{align}

Concerning the operator $\wt{\caln}^u$,
we first see that
\begin{align} \label{Lu-bdd}
   \sup\limits_{u\in \calu_{ad}}  \|G^u  \psi\|_{ \calv'}
   \leq \|u\|_{L^\9(H;H,\nu)} \|B\| \|\psi\|_{ \calh}.
\end{align}
Actually, by H\"older's inequality,
\begin{align*}
    |_{\calv}(\vf, G^u \psi)_{\calv'}|
    = \bigg| \int_{H} \l<B u, Q^\frac 12 D\vf\r> \psi d \nu \bigg|
  \leq& ( \int_{H} |\l<B u, Q^\frac 12 D\vf\r>|^2 d \nu )^\frac 12
         ( \int_{H} |\psi|^2 d \nu )^\frac 12 \\
  \leq& \|u\|_{L^\9(H;H,\nu)} \|B\| \|\vf\|_{\calv} \|\psi\|_{\calh},
\end{align*}
which yields \eqref{Lu-bdd}, as claimed.

Then,
arguing as above
and using \eqref{Integ-part} and \eqref{Lu-bdd}
we have that, for some $C>0$,
for any $\vf   \in \calv$,
\begin{align} \label{cala-calv}
   \|\wt{\caln}^{u} \vf \|_{\calv'} \leq C \|\vf\|_{\calv},
\end{align}
and
\begin{align} \label{cala-vf-vf}
 _{\calv}(\vf, \wt{\caln}^{u} \vf)_{\calv'}
  \leq&  -\frac 14 \|\vf\|_\calv^2.
\end{align}

Thus,
putting together the estimates above
we obtain \eqref{calna-v-v'} and \eqref{calna-dissp}.
\hfill $\square$

As a consequence of Proposition \ref{Prop-caln-cala}
and a classical result due to J.L. Lions
(see \cite{L69} or \cite[Theorem 4.10]{B10}),
we obtain that
there exist unique solutions $\wt\vf$ and $\wt p$ to \eqref{equa-wtvf}
and \eqref{equa-wtp}, respectively,
and so do the equations \eqref{equa-Nu} and \eqref{equa-p}.

\begin{theorem} \label{Thm-KE-Back}
Under the condition of Lemma \ref{Lem-integ-symm}.
Let $u\in \calu_{ad}$.
Then,  we have

$(i)$ For any objective function $g\in L^2(H,\nu)$,
there exists a unique solution $\vf^u$ to the Kolmogorov equation \eqref{equa-Nu}
such that
$\vf^u\in C([0,T]; \calh) \cap L^2(0,T; \calv)$,
$\frac{d}{dt} \vf^u \in L^2(0,T; \calv')$, and
\begin{align} \label{var-ke}
   \frac{d\vf^u}{dt}(t) = N_2 \vf^u(t) + \l<Q^\frac 12 Bu, D\vf^u\r>(t),
   \ \ a.e.\ t\in (0,T),\ \vf^u(0)=g,
\end{align}
where $\frac{d}{dt}$ is taken in the strong topology of $\calv'$,
or equivalently in $\mathscr{D}'(0,T; \calv')$.

$(ii)$ There exists a unique solution $p^u$ to
the backward equation \eqref{equa-p}
such that
$p^u\in C([0,T]; \calh) \cap L^2(0,T; \calv)$,
$\frac{d}{dt} p^u \in L^2(0,T; \calv')$,
and
\begin{align} \label{var-back}
    \frac{dp^u}{dt}(t) = -N_2 p^u(t) - G^u p^u(t) -1,
    \ \ a.e.\ t\in (0,T),\ p^u(T)=0.
\end{align}
\end{theorem}

The following result contains the uniform estimates and the continuity dependence
on controllers of the solutions to Kolmogorov equations.

\begin{theorem} \label{Thm-WP-KE}
Consider the situations as in Lemma \ref{Lem-integ-symm}.
We have

$(i)$
For any two  solutions  $\vf_1, \vf_2$ to \eqref{equa-Nu}
corresponding to the initial data $g_1$ and $g_2$, respectively,
we have
\begin{align} \label{vf12-g12}
   \|\vf_1-\vf_2\|_{C([0,T]; \calh)}
   + \|\vf_1- \vf_2\|_{L^2(0,T; \calv)}
   \leq C(\rho, T) \|g_1-g_2\|_\calh,
\end{align}
where $C(\rho, T)$ is independent of $u\in \calu_{ad}$.
In particular, one has
\begin{align} \label{unibdd-contr-sol}
   \sup\limits_{u\in \calu_{ad}}
   \|\vf^u\|^2_{C([0,T];\calh)}
   +  \sup\limits_{u\in \calu_{ad}}
     \int_0^T \int_{H} |Q^\frac 12 D\vf^u|_H^2 d \nu dt
     \leq C(\rho, T) <\9.
\end{align}

$(ii)$ For any $u, \wt{u} \in \calu_{ad}$
and $\lbb \in [0,1]$,
set $v :=  \wt{u} - u$.
Then, as $\lbb \to 0$,
\begin{align} \label{vflbb-vf-0}
   \|\vf^{u+\lbb v} - \vf^u\|_{C([0,T];\calh)}
   + \int_0^T \int_{H} |Q^\frac 12 D(\vf^{u+\lbb v} - \vf^u)|_H^2 d \nu dt  \to 0.
\end{align}
\end{theorem}

{\bf Proof. }
$(i)$
The estimate \eqref{vf12-g12}
follows from \eqref{equa-Nu} and
similar arguments as in the proof of \eqref{bdd-eN2u}.

$(ii)$
We replace $\vf$ by
$\vf_\lbb:= e^{-(4\rho^2 \|B\|^2 +1)t} (\vf^{u+\lbb v} - \vf^u)$
in  \eqref{equa-Nu} to obtain
\begin{align*}
   \frac{d}{dt} \vf_\lbb = \caln^{u} \vf_\lbb + \lbb e^{-(4\rho^2 \|B\|^2 +1)t}
                           \l<B v, Q^\frac 12 D \vf^{u+\lbb v}\r>,
\end{align*}
with $\vf_\lbb (0) =0$.
This, via \eqref{calna-dissp}, yields that
\begin{align*}
   \frac 12 \frac{d}{dt} \|\vf_\lbb\|_{\calh}^2
   \leq& -\frac 14 \|\vf_\lbb\|_\calv^2
         + \lbb e^{-(4\rho^2 \|B\|^2 +1) t} \int_{H} \l<B v, Q^\frac 12  D\vf^{u+\lbb v}\r> \vf_\lbb d \nu \\
   \leq&  -\frac 14  \|\vf_\lbb \|_\calv^2
         + 2 \lbb \rho \|B\|  \|\vf_\lbb\|_\calh(\int_{H} |Q^\frac 12 D\vf^{u+\lbb v}|_H^2 d \nu)^\frac 12.
\end{align*}
Thus, in view of the uniform boundedness \eqref{unibdd-contr-sol},
we obtain for some positive constant
$C'(T,\rho)$ independent of $\lbb$,
\begin{align*}
   \|\vf_\lbb\|^2_{C([0,T];\calh)}
   +  \int_0^T \int_{H} |Q^\frac 12 D\vf_\lbb|_H^2 d \nu dt
   \leq C'(T,\rho)\lbb^2 \to 0,\ as\ \lbb \to 0,
\end{align*}
which implies \eqref{vflbb-vf-0},
thereby finishing the proof.
\hfill $\square$\\

{\bf Proof of Theorem \ref{Thm-Contr-Symm}.}
Let $u_n$, $u_*$, $\vf_n$, $\vf_*$ and $\eta$
be as in the proof of Theorem \ref{Thm-Contr}.
Using \eqref{calna-v-v'} and \eqref{unibdd-contr-sol}
we obtain that
along a subsequence $\{n\}$,
\begin{align*}
    u_n \to u_*,\ \     &weak-star\ in\ L^\9(H; \nu),  \\
    \vf_n \to \vf_*,\ \ &weak-star\ in\ L^\9(0,T; \calh),   \\
                        &weakly\ in\ L^2(0,T;\calv),  \\
    \caln^{u_n} \vf_n\to \eta,\ \ &weakly\ in\ L^2(0,T; \calv').
\end{align*}
Then, using similar arguments as in the proof of Theorem \ref{Thm-Contr-Symm}
we can pass to the limits
in the dual pair $_{\calv_t}(\cdot, \cdot)_{\calv_t'}$
instead of the inner product $(\cdot, \cdot)_{\calh_t}$,
where $\calv_t := L^2(0,t; \calv)$ and
$\calv'_t$ is the dual space of $\calv_t$.

Hence, similarly to \eqref{0-eta-caln-rd}, we have that
\begin{align} \label{0-eta-calnu-symm}
    0\geq \int_0^T h(t)
          {}_{\calv'_t}(\eta - \caln^{u_*} \psi, \vf_*- \psi)_{\calv_t} dt
\end{align}
for any positive functions $h\in L^\9(0,T)$
and any $\psi \in \mathcal{M}$,
where $\mathcal{M}$ is defined as in the proof of Theorem \ref{Thm-Contr}.
Then,
since $\mathcal{M}$ is dense in $L^2(0,T;\calv)$
we can extend \eqref{0-eta-calnu-symm}
to all $\psi \in L^2(0,T; \calv)$.

Hence, taking $\psi = \vf_*- \ve f v$,
$f\in L^\9(0,T)$, $v\in \calf C_b^2(H)$,
similarly to \eqref{0-eta-calnu}
we obtain
\begin{align}
    0\geq \int_0^T h(t)
          {}_{\calv'_t}(\eta - \caln^{u_*} \vf_*, fv)_{\calv_t} dt
\end{align}
for any $f\in L^\9(0,T)$ and any $v\in \calf C_b^2(H)$,
which suffices to yield that
$\eta = \caln^{u_*} \vf_*$, $dt\times \nu$-a.e.,
thereby yielding that $u_*$ is an optimal controller for
Problem $(P^*)$.
Therefore, the proof is complete. \hfill $\square$ \\

{\bf Proof of Theorem \ref{Thm-Char-Optimal-Symm}.}
For any $u\in \calu_{ad}$, let $v : =u-u_*$
and $\vf^{u_*+\lbb v}$ be the solution to \eqref{equa-Nu} related to $u_*+\lbb v$,
$\lbb \in (0,1)$.

We infer from \eqref{equa-Nu} and \eqref{equa-p} that
\begin{align*}
   \frac{d}{dt} (\vf^{u_*+\lbb v} - \vf_*, p_*)
   = \lbb ( \l< B  v, Q^\frac 12D\vf^{u_*+\lbb v}\r>, p_*)
     - \int_{H} (\vf^{u_*+\lbb v} - \vf_*) d \nu,
\end{align*}
where $(\cdot, \cdot)$ denotes the inner product in $\calh$.
This yields that for any $\lbb \in (0,1)$,
\begin{align*}
   \int_0^T \int_{H}  \l<B(u-u_*), Q^\frac 12 D \vf^{u_*+\lbb v}\r> p_* d \nu dt
   =\frac 1 \lbb \int_0^T \int_{H} (\vf^{u_*+\lbb v} - \vf_*) d \nu dt
   \geq 0,
\end{align*}
where the last inequality is due to the optimality of $u_*$.

Therefore, taking $\lbb \to 0 $ and
using \eqref{vflbb-vf-0}, Fubini's theorem we obtain that
\begin{align*}
   \int_H \l<B(u-u_*), \int_0^T Q^\frac 12  D \vf_* p_* dt \r> d\nu \geq 0, \ \ \forall u \in \calu_{ad},
\end{align*}
which yields \eqref{u*-character},
thereby finishing the proof.
\hfill $\square$

\section{Optimal control of stochastic equations} \label{Sec-OP-SDE}

\subsection{General case}  \label{Subsec-OP-SDE-General}

In this subsection,
we first prove Theorem \ref{Thm-GWP-Mart} under Hypothesis $(H1)'$.
Then, we prove the first assertion $(i)$ of Theorem \ref{Thm-Contr-SDE}.
At last,
we show that Hypothesis $(H1)'$ can be implied by
the m-dissipativity of Kolmogorov operators,
i.e., Hypothesis $(H1)$-$(iii)$,
by applying the theory of generalized Dirichlet forms.  \\

{\bf Proof of Theorem \ref{Thm-GWP-Mart}.}
({\bf Existence})
By Hypothesis $(H1)'$,
there exists a conservative  Markov process
$\bbm=(\Omega, \mathscr{F}$, $(\mathscr{F}_t)_{t\geq 0}, (X(t))_{t\geq 0},$ $(\bbp_x)_{x\in H_0})$
such that $X(0)=x$, $\bbp$-a.s.,
the sample paths $t\mapsto \l<X(t), e_i\r>$ are continuous for every $i\geq 1$,
and for $\nu$-a.e. $x \in H$,
\begin{align} \label{Mart-N2}
    \vf(X(t)) - \int_0^t N_2 \vf(X(s)) ds
\end{align}
is an $(\mathscr{F}_t)$-martingale under $\bbp_x$
for all $\vf \in \mathcal{F} C_b^2(H)$.
Its transition semigroup $P_t$, $t>0$,
is given by $e^{tN_2}$, $t>0$
(see Remark \ref{Rem-i-ii-iii} $(iii)$),
where $N_2$ denotes the corresponding generator.

Let $\{e_i\} \subseteq D(\wt A)$ be the orthonormal basis of $H$
such that $Qe_i = q_i e_i$
with $q_i=0$ for $i\in J$
and $q_i >0$ for $i\notin J$,
where $J$ is a set of finitely many indices.
Set  $X_i(t) := \<X(t), e_i\>$,
$b_i(X(t)) := \<A (X(t)), e_i\>$
and
$(Bu)_i = \<Bu,e_i\>$,  $i \geq 1$.

For  every $i \notin J$, we set
\begin{align} \label{def-beta}
   \beta_i(t): = q_i^{-\frac 12} (X_i(t) - \int_0^t b_i(X(s))ds),\ \  t\geq 0.
\end{align}
By Definition \ref{Def-equacontr},
the sample paths $t\mapsto \beta_i(t)$ are continuous
for every $i\geq 1$ under $\bbp_x$ for $\nu$-a.e. $x\in H$.

Now, using standard regularization arguments
we infer from \eqref{Mart-N2} that
\begin{align} \label{Xi-J}
   X_i(t) = x_i + \int_0^t b_i(X(s)) ds, \ \ t\in [0,T], \ i\in J,
\end{align}
while $\beta_i$, $i \notin J$,
are independent  $(\mathscr{F}_t)$-Brownian motions
with $\beta_i(0) = q_i^{-\frac 12} x_i$ under $\bbp_x$, $x\in \ol{H}$.
(See, e.g.,  the proof of \cite[Corollary 1.10]{DRW09}.
Note that, unlike in \cite{DRW09},
the definition of $\beta_i$ in \eqref{def-beta} above is independent of $x$.)
Then, set
\begin{align*}
   M^u(t) := \exp(\sum\limits_{i\notin J}\int_0^t (B u)_i(X(s))d\beta_i(s)
                 - \frac 12 \sum\limits_{i\notin J} \int_0^t |(Bu)_i(X(s))|^2 ds),
\end{align*}
where $t\in [0,T]$.
Since ${\rm  sup}_{x\in H} |B u(x)|_H \leq \rho \|B\|$,
$\{M^u(t)\}$ is  an $(\mathscr{F}_t)$-martingale under $\bbp_x$
satisfying $\bbe M^u(T)=1$.
Hence, we have a new probability measure
$$\bbq^u_x :=   M^u(T) \cdot \bbp_x.$$
Then, Girsanov's theorem yields that
$$\wt{\beta}_{i}(t) := \beta_i(t) - \int_0^t ( B u)_i(X(s)) ds,\ \ t\in [0,T],\ i \notin J, $$
are independent $(\mathscr{F}_t)$-Brownian motions
with $\wt{\beta}_i(0) = q_i^{-\frac 12} x_i$
under $\bbq^u_x$.
Taking into account \eqref{def-beta}
we obtain  that
\begin{align} \label{equa-xi-perb}
   X_i(t) = \int_0^t b_i(X(s)) + q_i^\frac 12 ( B u)_i (X(s)) ds
   +   q_i^\frac 12 \wt{\beta}_{i}(t), \ \ t\in [0,T],\ i \notin J.
\end{align}
This together with \eqref{Xi-J} yields that $X(0)=x$, $\bbq^u_x$-a.s..
It also implies, via It\^o's formula, that under $\bbq^u_x$,
for any $\vf\in \calf C_b^2(H)$,
\begin{align}  \label{mart-fCb2}
 \vf(X(t)) - \int_0^t N_0 \vf(X(s)) + \l< B u(X(s)), Q^\frac 12  D \vf(X(s))\r>  ds
\end{align}
is an $(\mathscr{F}_t)$-martingale.

Hence,
$\bbm^u=(\Omega, \mathscr{F}$, $(\mathscr{F}_t)_{t\geq 0}, (X(t))_{t\geq 0},$ $(\bbq^u_x)_{x\in H_0})$
satisfies the property $(ii)$ of Definition \ref{Def-equacontr}.

It is also clear that $\bbm^u$ satisfies the property $(i)$ of Definition \ref{Def-equacontr},
since $P^u_t$ is bounded on $L^2(H,\nu)$.

Moreover,
we also have the Markov property for $(X(t))$ under $\bbq^u_x$,
i.e., for any $0<s,t<\9$,
\begin{align} \label{Mar-Qu}
   \bbe_{\bbq_x}(f(X(t+s))|\mathscr{F}_s)
  = \bbe_{\bbq_{X(s)}} (f(X(t))),\ \ \forall f\in \mathcal{B}_b(H).
\end{align}
To this end,
we first see that for any $F\in \mathcal{B}_b(H)$,
\begin{align*}
   \bbe_{\bbq_x}(F| \mathscr{F}_s)
   = \frac{\bbe_{\bbp_x} (FM^u(T)|\mathscr{F}_s)}{\bbe_{\bbp_x}(M^u(T)|\mathscr{F}_s)},\ \ 0\leq s\leq T.
\end{align*}
Then,
since $(M^u(t))$ is an $(\mathscr{F}_t)$-martingale under $\bbp_x$,
we have
\begin{align*}
   \bbe_{\bbq_x}&(f(X(t+s))|\mathscr{F}_s)
   = (M^u(s))^{-1}  \bbe_{\bbp_x}(f(X(t+s)) M^u(t+s)|\mathscr{F}_s) \\
   =&   \bbe_{\bbp_x}(f(X(t+s)) \exp(\sum\limits_{i\notin J}\int_s^{t+s} (B u)_i(X(r))d\beta_i(r)
                 - \frac 12 \sum\limits_{i\notin J} \int_s^{t+s} |(Bu)_i(X(r))|^2 dr) |\mathscr{F}_s),
\end{align*}
which along with the Markov property of $(X(t))$ under $\bbp_x$ yields that
\begin{align*}
   \bbe_{\bbq_x}(f(X(t+s))|\mathscr{F}_s)
   = \bbe_{\bbp_{X(s)}}(f(X(t)) M^u(t))
   =   \bbe_{\bbq_{X(s)}}(f(X(t))),
\end{align*}
where the last step is due to the martingale
property of $(M^u(t))$ under $\bbe_{\bbp_{X(s)}}$.
Thus, we obtain \eqref{Mar-Qu}, as claimed.

Therefore, we conclude that $\bbm^u$ is a solution to the martingale problem for \eqref{equa-X}
in the sense of Definition \ref{Def-equacontr}.

({\bf Uniqueness})
We adapt the arguments  as in the proof of \cite[ Theorem 8.3]{DR02}.
Let $\bbm'=(\Omega', \mathscr{F}'$, $(\mathscr{F}'_t)_{t\geq 0}, (X'(t))_{t\geq 0},$ $(\bbp'_x)_{x\in H_0})$
be another solution to the martingale problem of \eqref{equa-X},
with $(P^u_t)'$ and $(N_2^{u})'$ being the corresponding semigroup and generator, respectively.
Similarly,
let $(P^u_t)$ and $N^u_2$ be the semigroup and generator corresponding to $\bbm^u$.

We shall prove that for $\nu$-a.e. $x\in H$,
\begin{align} \label{P't-Pt}
    (P^u_t)' f(x) = P^u_t f(x) , \ \ t>0,\ \forall f\in C_b(H).
\end{align}

For this purpose,
we note that the property $(ii)$ in Definition \ref{Def-equacontr}
implies that under $\bbp_\nu:=\int_{H_0} \bbp_x \nu(dx)$,
$\vf(X'(t)) - \int_0^t N_0^u \vf(X'(s))  ds$
is a martingale for any $\vf \in \calf C_b^2(H)$.
It follows that,  for any $g\in L^2(H,\nu)$,
\begin{align*}
     &\int_{H_0} g(x) ((P^u_t)' \vf(x) - \vf(x) - \int_0^t (P^u_s)' N_0^u \vf(x) ds ) \nu(dx)  \\
   = & \bbe_v (g(X'(0))(\vf(X'(t)) - \vf(X'(0)) - \int_0^t N_0^u \vf (X'(s)) ds)) =0,
\end{align*}
which  implies that for any $\vf \in \calf C_b^2(H)$,
\begin{align*}
    (P^u_t)' \vf - \vf - \int_0^t (P^u_s)' N_0^u \vf ds =0,\ \ in\ L^2(H,\nu).
\end{align*}
Hence,
$\vf \in D((N_2^u)')$
and $ (N_2^u)' \vf = N^u_0 \vf$.

Taking into account $\calf C^2_b(H)$ is a core of $N_2^u$,
we obtain  that
$D(N^u_2) \subseteq D((N^u_2)')$ and $(N^u_2)' = N^u_2$ on $D(N_2)$.
But,
by Theorem \ref{Thm-Max-Nu},
$N_2^u$ is m-dissipative on $L^2(H,\nu)$.
Thus, we obtain $(N^u_2)' = N^u_2$,
which implies \eqref{P't-Pt}
and finishes the proof.
\hfill $\square$ \\

{\bf Proof of Theorem \ref{Thm-Contr-SDE} $(i)$.}
For any $u\in \calu_{ad}$,
let $X^u$ solve the martingale problem for \eqref{equa-X}
and $P_t^u$ be the corresponding transition semigroup,
i.e.,
$P_t^u f(x)= \bbe_{\bbp_x} f(X^u(t))$, $f\in L^2(H,\nu)$.

Then, for any $g\in D(N_2)$,
since $\calf C_b^2(H)$ is dense in $L^2(H,\nu)$,
we infer from Remark \ref{Rem-WP-Mart} $(ii)$ that
for any $t>0$,
$P^u_t g= e^{tN_2^u} g$, $\nu$-a.e. $x\in H$.
This yields that
\begin{align} \label{equi-bpxg-etnu}
   \int_0^T \int_H \bbe_{\bbp_x} g(X^u(t)) \nu(dx) dt
   = \int_0^T \int_H P^u_t g d\nu dt
   = \int_0^T \int_H e^{tN_2^u} g d\nu dt .
\end{align}

Thus, taking into account $\{e^{tN_2^u} g\}$ solves equation \eqref{equa-Nu}
in the space $L^2(H,\nu)$,
we infer that
the optimal controllers to Problem $(P^*)$
are also optimal to Problem $(P)$.

Actually,
let $u_*$ be an optimal controller to Problem $(P^*)$.
Suppose that  $u_*$ is not an optimal controller to Problem $(P)$,
then there exists $\wt u\in \calu_{ad}$ such that
\begin{align*}
   \int_0^T \int_H \bbe_{\bbp_x} g(X^{\wt u}(t)) \nu(dx) dt
   < \int_0^T \int_H \bbe_{\bbp_x} g(X^{u_*}(t)) \nu(dx) dt,
\end{align*}
which along with \eqref{equi-bpxg-etnu} yields that,
if $I_*$ denotes the minimum of objective functionals in Problem $(P^*)$,
\begin{align*}
  I_* \leq  \int_0^T \int_H e^{tN_2^{\wt u}} g d\nu dt
      < \int_0^T \int_H e^{tN_2^{u_*}} g d\nu dt
      =I_*,
\end{align*}
thereby yielding a contradiction.

Therefore, the proof is complete.
\hfill $\square$ \\

To end this subsection,
we show that the well-posedness of martingale problems
can be implied by the m-dissipativity of Kolmogorov operators in certain situations.

\begin{theorem} \label{Thm-H1iv}
Assume $(H1)$.
Assume additionally that $Tr Q<\9$ and $\int_H |A(x)|_H^2  \nu(dx)<\9$.
Then, $(H1)'$ holds, namely,
the martingale problem is well posed for \eqref{equa-X} when $u\equiv 0$.
\end{theorem}

In order to prove Theorem \ref{Thm-H1iv},
we  construct a nice Markov process
using the framework of \cite{BBR06},
which extends the generalized Dirichlet form in \cite{S99}
to $L^p$ spaces, $p\geq 1$.

We first see that, by $(H1)$ $(iii)$, since $L^2(H,\nu) \subseteq L^1(H,\nu)$,
$(N_0, \calf C_b^2(H))$
is also essentially m-dissipative in $L^1(H,\nu)$.

Then, we denote by $(N_1, D(N_1))$ and $G_\lbb^{(1)}:= (\lbb - N_1)^{-1}$
the closure of $(N_0, \calf C_b^2(H))$ in $L^1(H,\nu)$
and the corresponding resolvent, respectively, $\lbb >0$.
We say that $f$ is a 1-excessive function if $f\geq 0$
and $\lbb G_{1+\lbb}^{(1)} f \leq f$ for all $\lbb >0$.

\begin{lemma} \label{Lem-exces}
Consider the situations as in Theorem \ref{Thm-H1iv}.
Then,

$(i)$ For any $x\in H$,
\begin{align} \label{x2-G11}
     |x|_H^2 \leq G_1^{(1)} (|x|_H^2 + Tr Q + 2|A(x)|_H |x|_H) =: g(x).
\end{align}

$(ii)$ The function $g^{\frac 12}$ is 1-excessive.
\end{lemma}

{\bf Proof.}
$(i)$
Define the projection operator $P_n$ by
$P_nx := \sum_{i=1}^n \l<x,e_i\r> e_i$, $x\in H$.
Using similar regularization procedure as in the proof of \cite[Lemma 5.5]{BBDR06},
we see that $|P_n x|_H^2 \in D(N_1)$,
and
\begin{align*}
   (1-N_1) |P_n x|_H^2 = |P_n x|_H^2 - \sum\limits_{i=1}^n q_i - 2\<A(x), P_n x\>.
\end{align*}
Since $Tr Q <\9$
and $\int_H |A(x)|_H^2 + |x|_H^2 \nu(dx)<\9$,
we obtain
\begin{align*}
     |P_n x|_H^2 \leq G_1^{(1)} ( |P_n x|_H^2 + Tr Q + 2|A(x)|_H |P_n x|_H) \in L^1(H,\nu),
\end{align*}
which implies \eqref{x2-G11} by passing to the limit.

$(ii)$ By the resolvent equation we have $\lbb G_{1+\lbb}^{(1)} g \leq g$.
Then, using Jensen's inequality we obtain
\begin{align*}
   \lbb G_{1+\lbb}^{(1)} g^\frac 12
   \leq \frac{\lbb}{\lbb +1} ((\lbb +1) G_{1+\lbb}^{(1)} g)^\frac 12
   = \frac{\lbb^\frac 12}{(\lbb+1)^\frac 12} (\lbb G_{1+\lbb}^{(1)} g)^\frac 12
   \leq g^\frac 12,
\end{align*}
which finishes the proof.
\hfill $\square$\\

{\bf Proof of Theorem \ref{Thm-H1iv}.}
The proof is based on \cite[Theorem 1.1]{BBR06}.
We first see that
the condition $(II)$ of \cite[Theorem 1.1]{BBR06} is satisfied
with $\cala = \calf C_b^2(H)$.

Moreover,
$F_n:= \{x\in H, |x|_H \leq n\}$ is weakly compact  in $H$, $n\geq 1$,
and
\begin{align*}
   R_1 (I_{F_n^c}) \leq \frac 1n \int_H g^\frac 12 d\nu \to 0,\ \ as\ n\to \9,
\end{align*}
where $R_1$ is  defined as in \cite{BBR06},
$I_{F_n^c}$ denotes the characteristic function of the complement set of $F_n$,
and $g$ is the 1-excessive function as in Lemma \ref{Lem-exces}.
Then, using \cite[Remark 2.2]{BBR06}
(with $f_0 =1$, $V_\beta = G_1^{(1)}$, $\beta=1$)
we have that $\{F_n\}_{n\geq 1}$ is a $\nu$-nest of weakly compact sets,
and so the condition $(I)$ of \cite[Theorem 1.1]{BBR06} is also satisfied.

Thus,
by virtue of \cite[Theorem 1.1]{BBR06},
we obtain a $\nu$-standard right process
$\bbm = (\Omega, \mathscr{F},(\mathscr{F}_t)_{t\geq 0}, (X(t))_{t\geq 0}, (\bbp_x)_{x\in H})$
in the state space $H$
with c\`{a}dl\`{a}g sample path in the weak topology  of $H$.
(Note that,
the life time of $(X(t))$ is infinite since $N_1 1=0$.)

To show that in our case the paths
$t\mapsto \l<X(t), e_i\r>$, $i\geq 1$,
are continuous
we adapt a method developed in \cite[Section 6]{DR02}.
So, let $\{\vf_n, n\in \bbn\}$
be a countable subset of $\calf C_b^2(H)$ separating the points of $H$.
Then, by exactly the same arguments as in the proof of \cite[Theorem 6.3]{DR02}
one obtains that for all $n\in \bbn$, $s<t$,
\begin{align*}
    \int\limits_{\Omega} |\vf_n(X(t)) - \vf_n(X(s))|^4 d\nu
    \leq C_n (t-s)^{\frac 32},
\end{align*}
where $C_n \in (0,\9)$.
Since we already know that $X(t)$, $t\geq 0$,
is  weakly c\`{a}dl\`{a}g $\bbp_x$-a.s. for $\nu$-a.e. $x\in H$,
this together with the proof of Kolmogorov's continuity criterion implies that
\begin{align} \label{Lamda-P-1}
    \bbp_\nu (\Lambda_0) =1,
\end{align}
and so
\begin{align}
   \bbp_x(\Lambda_0) =1, \ \ for\ \nu-a.e.\ x,
\end{align}
where $\Lambda_0:= \bigcap_{k,l\in \bbn} A_k^{(l)}$ with
\begin{align}
     A^{(l)}_k
     :=& \big\{w\in \Omega: \exists n_0 \in \bbn, \forall n\geq n_0, \forall s,t\in D_n \cap [0,l],
         |s-t| \leq 2^{-n_0}: \nonumber \\
       & \qquad \qquad  |\vf_n(X(t)) - \vf_n(X(s))| \leq 2^{-k}  \big\}, \nonumber  \\
   D_n:=& \{k 2^{-n}, k\in \bbn \cup \{0\}\},\ \ D:= \bigcup_{n\in \bbn} D_n.\nonumber
\end{align}
This yields that for $\nu-a.e.\ x$,
under $\bbp_x$
the paths $t\mapsto \vf_n(X(t))$ are continuous  for all $n\geq 1$,
and so are the paths $t\mapsto \l<X(t), e_i\r>$ by density, $i\geq 1$.

Therefore,
we conclude that $\bbm$ solves the martingale problem of \eqref{equa-X}
in the sense of Definition \ref{Def-equacontr}.

The uniqueness can be proved by using similar arguments as in the proof of Theorem \ref{Thm-GWP-Mart}.
The proof is complete.
\hfill $\square$

\subsection{Symmetric case} \label{Subsec-OP-SDE-Symm}

In this case,
the nice feature is that
the associated Dirichlet forms are coercive closed forms.
This enables us to apply the general framework of Dirichlet forms as in \cite{MR92}
to obtain that,
the corresponding semigroup is even holomorphic
and Hypothesis $(H1)'$ also holds,
i.e., the martingale problem is well posed for \eqref{equa-X} when $u\equiv 0$.

Below we fix $\lbb > 4\rho^2 \|B\|^2$.
For any $u\in \calu_{ad}$,
we define the bilinear map $\cale_\lbb^u: \calf C_b^2(H) \times \calf C_b^2(H) \to \bbr$ as follows
\begin{align*}
   \cale_\lbb^u(\vf, \psi):= \frac 12 \int_H \l<Q^\frac 12 D \vf, Q^\frac 12 D \psi\r> d \nu
 - ( \l<Bu, Q^\frac 12 D\vf\r>, \psi)
 + \lbb (\vf, \psi)
\end{align*}
for any $\vf, \psi \in \calf C_b^2(H)$,
where $(\cdot,\cdot)$ is the inner product in $\calh:= L^2(H,\nu)$.
Under the closability assumption of $(Q^\frac 12 D, \calf C_b^1(H))$,
we can extend $(\cale_\lbb^u, \calf C_b^2(H))$
to the  closed form $(\cale_\lbb^u, \calv)$,
where $\calv := W^{1,2}(H,\nu)$.

\begin{lemma} \label{Lem-DF-symm}
Assume $(H1)$.
Assume that $N_2$ is symmetric
and $(Q^\frac 12D, \mathcal{F}C_b^1(H))$
is closable from $L^2(H,\nu)$ to $L^2(H;H,\nu)$.
Then, $(\cale_\lbb^u, \calv)$ is a coercive closed form.
\end{lemma}

{\bf Proof.}
We need only to check that $(\cale_\lbb^u, \calv)$
satisfies the weak sector condition,
namely,  for some $K>0$,
\begin{align} \label{veu-esti}
   \cale_{1+\lbb}^u(\vf, \psi)
   \leq K  \cale_{1+\lbb}^u(\vf, \vf)^\frac 12  \cale_{1+\lbb}^u(\psi, \psi)^\frac 12,\ \ \forall \vf, \psi \in \calv.
\end{align}
For this purpose,
it suffices to prove that
for some $c>0$
\begin{align} \label{cale-lowbbd}
   \cale_\lbb^u (\vf, \vf) \geq c \|\vf\|_\calv^2,\ \ \forall \vf \in \calv.
\end{align}
In order to prove \eqref{cale-lowbbd},
since $\lbb> 4 \rho^2 \|B\|^2$,
using Cauchy's inequality we get
\begin{align*}
  |( \l<Bu, Q^\frac 12 D\vf\r>, \vf) |
  \leq& \rho \|B\| \|Q^\frac 12 D\vf\|_{L^2(H; H,\nu)} \|\vf\|_{\calh} \\
  \leq& \frac 14 \|Q^\frac 12 D\vf\|^2_{L^2(H;H,\nu)}
         + 4 \rho^2 \|B\|^2 \|\vf\|_\calh^2,
\end{align*}
which implies that
\begin{align*}
   \cale_\lbb^u(\vf, \vf)
   \geq \frac 14 \|Q^\frac 12 D \vf\|^2_{L^2(H;H,\nu)}
        + (\lbb - 4 \rho^2 \|B\|^2) \|\vf\|_\calh^2,
\end{align*}
thereby yielding \eqref{cale-lowbbd}
with $c= \min\{\frac 14, (\lbb -2\rho^2 \|B\|^2)\}>0$, as claimed.
\hfill $\square$

Now,
by virtue of \cite[I. Proposition 2.16]{MR92},
we have
the one-to-one correspondence between $(\cale_\lbb^u, \calv)$
and the generator $(L^u_\lbb, D(L^u_\lbb))$,
where $L^u_\lbb$ is the unique element in $\calh$
such that $(-L^u_\lbb \vf, \psi) = \cale_\lbb^u(\vf, \psi)$
for all $\psi \in \calv$
and $\vf \in D(L^u_\lbb) := \{\vf\in \calv, \psi \to \cale_\lbb^u(\vf, \psi)\ is\ continuous\ w.r.t.\ (\cdot,\cdot)^\frac 12\ on\ \calv\}$.
Since $L^u_\lbb$ and $N^u_2 - \lbb $ coincides on $\calf C_b^2(H)$,
it follows that
$L^u_\lbb = N^u_2 - \lbb$ and $D(L^u_\lbb) = D(N_2)$.

The following result states that the corresponding semigroup is holomorphic,
which enables one to solve equation \eqref{equa-X}
in the space $\calh$
and also the optimal control problems
even for the objective functions in the space $\calh$.

\begin{corollary}  \label{Cor-RegPt-symm}
Consider  the situations as in Lemma \ref{Lem-DF-symm}.
Let $(e^{tN_2^u})$ be the semigroup corresponding to $(N_2^u, D(N_2))$.
Then, for all $t>0$
and for any $g\in \calh$,
we have $e^{tN_2^u}g\in D(N_2)$.
In particular,
$e^{tN_2^u}g$ is the unique solution to \eqref{equa-Nu}.
\end{corollary}

{\bf Proof.}
By virtue of
I. Corollary 2.21 and I. Theorem 2.20 of \cite{MR92},
we have that
$L^u_\lbb$ generates a holomorphic semigroup $(e^{t L^u_\lbb})$
on some sector in $\bbc$
such that for all $t>0$ and $g\in \calh$,
$e^{t L^u_\lbb} g\in D(L^u_\lbb)$
and so
$ e^{tN_2^u} g \in D(N_2)$,
due to
$e^{t L^u_\lbb} = e^{tN_2^u} e^{-\lbb t}$
and $D(L^u_\lbb) = D(N_2)$.
This yields that $e^{tN_2^u} g$ solves \eqref{equa-X} in the space $\calh$.
Moreover, the uniqueness of solutions to \eqref{equa-Nu} follows from the monotonicity
of $N_2^u$.
Therefore, the proof is complete.
\hfill $\square$

Below we show that $(H1)'$ can be implied from $(H1)$ in the symmetric case.

\begin{theorem} \label{Thm-H1'-H1-Sym}
Consider the situations as in Lemma \ref{Lem-DF-symm}.
Then, $(H1)'$ holds,
namely,
the martingale problem is well posed for \eqref{equa-X} when $u\equiv 0$.
\end{theorem}

{\bf Proof.}
We construct the Markov process by using the framework of Dirichlet forms in \cite{MR92}.
First we see that, when $u\equiv 0$, $\lbb =0$,
\begin{align*}
   \cale(\vf, \psi)
   (:=  \cale_0^0(\vf, \psi))
   = \frac 12 \sum\limits_{k=1}^\9 q_k \int_H \p_k \vf \p_k \psi d \nu,\ \ \vf, \psi \in \calf C_b^2(H).
\end{align*}
Since $Q$ is bounded on $\calh$,
$\sup_{i\geq 1} q_i <\9$.
Then,
taking into account $\calf C_b^2(H) \subseteq \calv$ is dense and separates points of $\calh$,
and using similar arguments
as in the proof of \cite[IV. Proposition 4.2]{MR92}
we have the quasi-regularity of $(\cale, \calv)$.

Hence,
by virtue of \cite[IV. Theorem 3.5]{MR92},
we obtain a $\nu$-tight special standard process $\bbm$
associated with $(\cale, \calv)$ hence also with $(N_2, D(N_2))$,
and its life time $\zeta =\9$ since $N_2 1 =0$.

Since the semigroup $e^{tN_2^u}$, $t\geq 0$,
is bounded on $\calh$,
the first property $(i)$ of Definition \ref{Def-equacontr} holds.

Moreover,
since $(\cale, \calv)$ has the local property (see \cite[V. Definition 1.1]{MR92}),
\cite[V. Theorem 1.5]{MR92} yields
that the sample path of $\bbm$ is continuous.
In view of \cite[Proposition 1.4]{BBR06}
and Remark \ref{Rem-A-bdd} $(i)$,
we also infer that
the property $(ii)$ in Definition \ref{Def-equacontr} is   satisfied for $\bbm$.

Thus, $\bbm$ solves the martingale problem for \eqref{equa-X}
when $u\equiv 0$.

The uniqueness of solutions to martingale problem can be proved similarly as
in Theorem \ref{Thm-GWP-Mart}.
Therefore,
the proof is complete.
\hfill $\square$ \\

{\bf Proof of Theorem \ref{Thm-Contr-SDE} $(ii)$}
For any $u\in \calu_{ad}$,
let $X^u$ solve the martingale problem for \eqref{equa-X}.
Similarly to \eqref{equi-bpxg-etnu},
for any $g\in \calh$ we have
\begin{align} \label{equi-bpxg-etnu-symm}
   \int_0^T \int_H \bbe_{\bbp_x} g(X^u(t)) d\nu dt
   = \int_0^T \int_H e^{tN_2^u} g d\nu dt.
\end{align}
Moreover, since by Corollary \ref{Cor-RegPt-symm} $e^{tN_2^u} g \in D(N_2)$, $t>0$,
we have
\begin{align} \label{equa-N2u-symm}
   \frac{d}{dt} e^{tN_2^u} g = N_2^u e^{tN_2^u} g,\ \ t\in (0,T],\ in\ \calh.
\end{align}

Let $\vf^u$ be the variational solution to \eqref{equa-X} as in Theorem \ref{Thm-KE-Back}.
Then, since $D(N_2) \subseteq \calv$,
using Lemma \ref{Lem-integ-symm}
and arguing as in the proof of \eqref{bdd-eN2u} we get
\begin{align*}
   \frac 12 \frac{d}{dt} \|e^{tN_2^u} g - \vf^u(t)\|^2_\calh
   =  _{\calv}(e^{tN_2^u} g - \vf^u(t), N_2^u(e^{tN_2^u} g - \vf^u(t)))_{\calv'}
   \leq C\|e^{tN_2^u} g - \vf^u(t)\|^2_\calh,
\end{align*}
which, via Gronwall's inequality, yields
$\|e^{tN_2^u} g - \vf^u(t)\|^2_\calh =0$ for any $t\in [0,T]$
and so,
\begin{align}
   e^{tN_2^u} g  = \vf^u(t), \ \ \nu-a.e.\ x,\ t\in[0,T].
\end{align}
This along with \eqref{equi-bpxg-etnu-symm} yields
\begin{align}
   \int_0^T \int_H \bbe_{\bbp_x} g(X^u(t)) d\nu dt
   = \int_0^T \int_H \vf^u(t) d\nu dt.
\end{align}

Therefore,
the optimal controllers to Problem $(P^*)$
are also optimal to Problem $(P)$.
The proof is complete.
\hfill $\square$

\section{Applications} \label{Sec-Application}

\subsection{Singular  dissipative stochastic equations} \label{Subsec-SingSDE}

We consider the singular dissipative stochastic equation as in \cite{DR02}
\begin{align} \label{equa-X-sing}
   &dX(t) = \wt{A}X(t)dt + F(X(t)) dt + Q^\frac 12 B u(X(t)) dt + Q^\frac 12  dW(t), \ \ t\in (0,T), \nonumber \\
   &X(0) = x \in H.
\end{align}
Here,
$\wt{A} : (D(\wt{A})) \subseteq H \mapsto H$ is
m-dissipative linear operator and
$F:D(F) \subset H \rightarrow H$ is an m-dissipative singular valued operator,
i.e.,
$\<F(x)-F(y), x-y\> \leq 0$, $\forall x,y\in D(F)$,
and
${\rm Range} (I-F):= \bigcup_{x\in D(F)} (x-F(x)) = H$.
The operators $B$ and $Q$ are as in \eqref{equa-X},
with the orthonormal basis $\{e_i\} \subseteq D(\wt A)$,
defined by $Qe_i = q_i e_i$, $q_i>0$, $i\geq 1$.

Let $A$ be defined by \eqref{def-A}.
We have
\begin{align*}
   \<A(x), D\vf(x)\> = \l<x, \wt A D\vf(x)\r> +\<F(x), D\vf(x)\>, \ \ \forall \vf \in \calf C_b^2(H),
\end{align*}
and $D(A) = D(F)$.

Let us recall the functional framework in \cite{DR02}.
Let $\mathcal{E}_{\wt{A}}(H)$ be the linear span of all (real parts of)
functions of the form $\vf = e^{i\<h,\cdot \>} $, $h\in D(\wt{A})$.

In addition,
we shall assume that
\begin{enumerate}
  \item[(A1)] $(i)$ There exists $\omega>0$ such that
  \begin{align*}
      \l<\wt{A} x, x\r> \leq - w |x|_H^2, \ \ \forall x\in H.
  \end{align*}
  $(ii)$ $Q$ is  bounded, self-adjoint and positive definite,
  $Q^{-1} \in L(H)$
  and  ${\rm Tr} \wt{Q}<\9$, where
  \begin{align*}
     \wt{Q}x: = \int_0^\9 e^{t\wt{A}}Qe^{t\wt{A}} x dt, \ \ x\in H.
  \end{align*}
    \item[(A2)] There exists a Borel probability measure $\nu$ on $H$ such that

  $(i)$ $\int_{D(F)} (|x|_H^{12} + |F(x)|_H^2 + |x|_H^4|F(x)|_H^2) \nu(dx) <\9. $

  $(ii)$ For all $\vf \in \mathcal{E}_{\wt{A}}(H)$ we have $N_0\vf \in L^2(H, \nu)$ and
  \begin{align*}
             \int_{H} N_0 \vf(x) \nu(dx) =0.
  \end{align*}

  $(iii)$ $\nu(D(F)) =1$.
\end{enumerate}

We see that Assumption $(A2)$
implies Hypothesis $(H1)$ $(i)$ and $(ii)$.
It also follows from \cite[Theorem 2.3]{DR02} that,
under Assumptions $(A1)$ and $(A2)$,
$(N_0, \mathcal{E}_{\wt{A}}(H))$ is essentially m-dissipative in $L^2(H, \nu)$,
and so is $(N_0, \mathcal{F}C_b^2 (H))$.
Thus Hypothesis $(H1)$  is satisfied.

Moreover,
\cite[Theorem 7.4]{DR02} yields that the
martingale problem is well posed for \eqref{equa-X}
in the case $u\equiv 0$,
which yields $(H1)'$.

Thus, both
Hypotheses $(H1)$ and $(H1)'$ are fulfilled.

As regards the closability of $D$ we have
\begin{proposition} \label{Prop-D-close}
Assume $(A1)$ and $(A2)$.
Then,
$D$ is closable from $L^2(H, \nu)$ to $L^2(H; H, \nu)$.
\end{proposition}

{\bf Proof.}
Let $\{\wt{e}_k\}_{k\geq 1}$ be an orthonorm basis of $H$
such that
$\wt{Q}\wt{e}_k= \wt{q}_k \wt{e}_k$, $\wt q_k>0$, $k\geq 1$,
and set $x_k := \<x,\wt{e}_k\>$,
$D_k \vf := \<D \vf, \wt{e}_k\>$, $x\in H$.

Let $\mu$ be the Gaussian measure with mean zero and covariance operator $\wt Q$.
We have the integration by parts formula,
\begin{align} \label{Integ-part-Dk-mu}
   \int_H D_k \vf \psi d\mu
   = - \int_H \vf D_k \psi d\mu
       + \frac{1}{\wt{q}_k} \int_H x_k \vf \psi d\mu, \ \ \forall \vf, \psi \in \calf C_b^2(H).
\end{align}

Moreover, for the infinitesimal invariant measure $\nu$,
it is known that
$\nu = \rho \cdot \mu$
with $\rho^\frac 12 \in W^{1,2}(H, \mu)$
(see \cite[page 292]{DR02}).
Note that $\rho \in \ol{\calf C_b^2(H)}^{W^{1,1}(H,\mu)}$,
$D_k(\rho \vf) = D_k \rho \vf + \rho D_k \vf$ for any $\vf \in \calf C_b^2(H)$,
$D\rho = 2 \rho^\frac 12 D(\rho^\frac 12)$ in $L^1(H; H,\nu)$,
and so $D_k \rho / \rho \in L^2(H, \nu)$.

Now, we take any $(\vf_n) \subseteq \calf C_b^2(H)$
such that
\begin{align} \label{vfn0-Dvfng}
    \vf_n \to 0,\ in\ L^2(H,\nu),\  \ D_k\vf_n \to g, \ \ in\ L^2(H, \nu).
\end{align}
We shall prove that
\begin{align} \label{g-0}
   g(x)=0,\ \ \nu-a.e.\ x.
\end{align}

For this purpose,
for any $\psi \in \calf C_b^2(H)$
we set $\psi_{\ve, k} := (1+\ve|x_k|^2)^{-1}\psi$, $\ve >0$.
Note that, by \eqref{vfn0-Dvfng},
\begin{align} \label{psivek-g-rho}
   \int_H \psi_{\ve, k} g \rho d\mu
   =& \lim\limits_{n\to \9} \int_H\psi_{\ve, k} D_k \vf_n \rho d\mu \\
   =&  \lim\limits_{n\to \9} \int_H D_k (\psi_{\ve, k}\vf_n \rho) d\mu
      - \lim\limits_{n\to \9} \int_H \vf_n  (D_k  \psi_{\ve, k}\rho  + \psi_{\ve, k}  D_k  \rho ) d\mu  \nonumber
\end{align}
Since $\psi_{\ve, k} \in \calf C_b^2(H)$
and $D_k\rho / \rho \in L^2(H,\nu)$,
using \eqref{vfn0-Dvfng} we see that the last limit on the right hand side above
equals to zero.
Regarding the remaining limit,
we take a sequence $\rho_m \in \calf C_b^2(H)$, $m\geq 1$,
such that
$\rho_m \to \rho$ in $W^{1,1}(H,\mu)$.
Then, using the integration by parts formula \eqref{Integ-part-Dk-mu}
we have
\begin{align*}
   \int_H D_k (\psi_{\ve, k} \vf_n \rho) d\mu
   =& \lim\limits_{m\to \9} \int_H D_k (\psi_{\ve, k} \vf_n \rho_m) d\mu \\
   =& -\frac{1}{\wt q_k} \lim\limits_{m\to \9} \int_H x_k \psi_{\ve, k} \vf_n \rho_m d\mu \\
   =& -\frac{1}{\wt q_k} \int_H x_k \psi_{\ve, k} \vf_n \rho  d\mu,\ \ as\ m\to \9,
\end{align*}
where in the last step we used the fact that
$\sup_{x\in H} |x_k \psi_{\ve,k}(x)| \leq C_\ve <\9$.
Hence, using again \eqref{vfn0-Dvfng}
and the boundedness of  $\sup_{x\in H} |x_k \psi_{\ve,k}(x)|$
we get
\begin{align*}
   \lim\limits_{n\to \9}
   \int_H D_k (\psi_{\ve, k} \vf_n \rho) d\mu
   = -\frac{1}{\wt q_k} \lim\limits_{n\to \9}
   \int_H x_k \psi_{\ve, k} \vf_n \rho d\mu  =0.
\end{align*}
Thus,
combing back to \eqref{psivek-g-rho}
we obtain
\begin{align*}
   \int_H \psi_{\ve, k} g \rho d\mu  =0.
\end{align*}
Taking the limit $\ve\to 0$ yields that,
for any $\psi \in \calf C_b^2(H)$,
\begin{align*}
   \int_H \psi  g \rho d\mu  = 0,
\end{align*}
which yields \eqref{g-0}
and finishes the proof.
\hfill $\square$

The compact embedding of $W^{1,2}(H,\nu)$ to $L^2(H,\nu)$ also holds in certain situations.
Following \cite{DRW09} we assume additionally that
\begin{enumerate}
  \item[(A3)] \ \ \ \
  $(i)$ $\int_0^\9 (1+t^{-\a}) \|e^{t\wt A} \sqrt{Q} \|^2_{HS} dt <\9$ for some $\a>0$,
  where $\|\cdot\|_{HS}$ denotes the norm on the space of all Hilbert-Schmidt operators on $H$,
  and $Q^{-\frac 12 } \in L(H)$.

  $(ii)$ $(1+w- \wt{A}, D(\wt{A}))$ satisfies the weak sector condition,
  i.e., for some $K>0$,  for any $x,y\in D(\wt A)$,
  \begin{align*}
     \l<(1+ w - \wt{A}) x, y\r>
     \leq K  \l<(1+ w - \wt{A}) x, x\r>^\frac 12   \l<(1+w  - \wt{A}) y, y\r>^\frac 12.
  \end{align*}

  $(iii)$ There exists a sequence of $\wt{A}$-invariant finite dimensional subspace
  $H_n \subseteq D(\wt{A})$ such that
  $\cup_{n=1}^\9 H_n$ is dense in $H$.
\end{enumerate}
It follows from \cite[Theorem $1.6$]{DRW09} that the
Harnack inequality holds for the $\nu$-version transition semigroup $p^\nu_t$
corresponding to \eqref{equa-X-sing} when $u\equiv 0$.
In particular, by \cite[Corollary 1.9]{DRW09},
$p^\nu_t$ has a density with respect to $\nu$
and is even hyperbounded,
i.e., $\|p^\nu_t\|_{L^2(H,\nu) \to L^4(H, \nu)} <\9$
for some $t>0$.

Thus, by virtue of \cite[Theorem 1.2]{GW02},
we obtain the compactness of $p^\nu_t$
and also the compact embedding of $W^{1,2}(H,\nu)$ into $L^2(H, \nu)$.
Hence, Hypothesis $(H3)$ is satisfied.

Now,
we conclude from Theorem \ref{Thm-Contr-SDE} in Section \ref{Sec-Main} that
\begin{theorem} \label{Thm-Contr-Sing}
Consider the controlled stochastic singular differential equation \eqref{equa-X-sing}.
Assume Hypotheses $(A1)$ and $(A2)$.
Assume additionally $(H2)$ or $(A3)$.
Then, for any $g\in D(N_2)$,
there exists an optimal control $u_*$ for the optimal control problem below
\begin{align*}
  {\rm Min} \bigg\{\int_0^T \int_H & \bbe_{\bbp_x}g(X^u(t)) \nu(dx) dt;\ u\in \calu_{ad},
                    \ \bbp_x\circ (X^u)^{-1}\ solves \nonumber \\
              & the\ martingale\ problem\ for\ \eqref{equa-X-sing}\ for\ \nu-a.e.\ x\in H \bigg\}.
\end{align*}
In particular,
under Hypotheses $(A1)$, $(A2)$ and $(A3)$,
we can take $B = I_d$.

\end{theorem}

As a specific example of \eqref{equa-X-sing},
we consider the controlled gradient system
\begin{align}
   d X  &= \wt{A} X dt + \p U(X) dt + Bu(X) +  dW(t), \label{equa-X-grad} \\
   X(0) &= x \in H.  \nonumber
\end{align}
Here, we take $Q = Id$,
$\wt A, B$ are the operators as in \eqref{equa-X-sing}
satisfying additionally that $\wt A^{-1}$ is of trace class.
and $\p U$ denotes the subdifferential of  a
convex and lower semicontinuous function $U: H \to (-\9, \9]$,
satisfying that
$\{U<\9\}$ is open, $\mu(\{U<\9\}) >0$
and
\begin{align*}
   \rho := Z^{-1} e^{-2U(x)} \in L^1(H, \mu),
\end{align*}
where $\mu$ is the Gaussian measure of mean zero and
covariance operator $-\frac 12 \wt A^{-1}$
and $Z:= \int_H e^{-2U(x)} \mu(dx)$.

We know from   \cite[ Section 9.1]{DR02} that
Assumptions $(A1)$ and $(A2)$ are fulfilled
and,
in particular,
the  Kolmogorov operator $N_2$ is symmetric.

Therefore,
by virtue of Theorem \ref{Thm-Contr-SDE},
for  more general objective functions $g\in L^2(H,\nu)$
we have
the existence as well as first-order necessary condition \eqref{u*-character}
of the feedback control problem $(P)$ for the gradient system \eqref{equa-X-grad}.

\subsection{Stochastic reaction-diffusion equation} \label{Subsec-RDE}

Consider the controlled stochastic reaction-diffusion equation below
as in \cite{DP04}
\begin{align} \label{equa-rd-X}
   dX =& \Delta X dt- p(X) dt + C^\frac 12 B u(X) dt + C^\frac 12 dW,   \\
   X(0)=& x \in H,    \nonumber
\end{align}
where
$H=L^2(\mathcal{O})$, $\mathcal{O} = [0,1]$,
$\Delta$ is the realization of the Laplace operator
with Dirichlet boundary condition, i.e.,
$D(\Delta) = H^2(\mathcal{O}) \cap H_0^1(\mathcal{O})$,
$B$ is a bounded operator on $H$,
and $W$ is a cylindrical Wiener process on $H$,
$W(t) = \sum_{k=1}^\9 e_k \beta_k(t)$
is a cylindrical Wiener process
on a stochastic basis $(\Omega, \mathscr{F}, (\mathscr{F}_t)_{t\geq 0}, \bbp)$,
where $e_k$ are the eigenbasis of $-\Delta$,
such that
$-\Delta e_k = \lbb_k e_k$,
$\lbb_k\geq 0$, $k\geq 1$.

Concerning $p$ and $C$ we assume that
\begin{enumerate}
  \item[(B)] \ \ \ \
  $(i)$ $p$ is a polynomial of degree $d>1$,
  its derivative $p'(\xi) \geq 0$, $\forall \xi \in \bbr$.

  $(ii)$ $C= (-\Delta)^{- \g}$, $\g> -\frac 12 $.
\end{enumerate}

In this case, we have $A(x) = \Delta x - p(x)$
and  $D(A) = \{x\in L^{2d}(\mathcal{O})\}$.

When $u\equiv 0$,
it is known (see Theorem 4.8 of \cite{DP04}) that,
for each $x\in H$,
there exists a unique generalized solution $X(\cdot,x)$ to \eqref{equa-rd-X}.

Moreover, by \cite[Theorem 4.16]{DP04}),
the transition semigroup $P_t: \mathcal{B}_b(H) \to \mathcal{B}_b(H)$ defined by
$(P_t \vf)(x) = \bbe \vf(X(t,x))$,  $x\in H$, $\vf\in C_b(H)$,
has a unique invariant measure $\nu$
satisfying that
\begin{align}
  \lim\limits_{t\to \9}
  P_t \vf (x) = \int_H \vf(y) \nu(dy)
\end{align}
and (see \cite[Proposition 4.20]{DP04})
\begin{align} \label{x2d}
   \int_H |x|_{L^{2d}(\mathcal{O})}^{2d}  \nu(dx) <\9.
\end{align}
Furthermore, from \cite[Section 4.6]{DP04}
we have that $P_t$ can be uniquely extended to a $C_0$-semigroup of contractions on $L^2(H;\nu)$.
By Theorem 4.23 of \cite{DP04},
the infinitesimal generator $N_2$ of $P_t$ is the closure in $L^2(H, \nu)$
of the operator
\begin{align} \label{N0-rde}
   (N_0 \vf)(x)
   := \frac 12 Tr [(-\Delta)^{-\g} D^2 \vf](x)
     + \<x, \Delta D \vf\>
     - \<p(x), D\vf\>,
\end{align}
where $x\in H$,
$\vf\in \mathcal{E}_\Delta(H)$ with
$\mathcal{E}_\Delta(H)$ defined similarly as
$\mathcal{E}_{\wt A}(H)$ in the previous subsection.

Now, let us check the Hypothesis $(H1)$.
We first infer from \eqref{x2d} that Hypothesis $(H1)$ $(i)$ is satisfied.
Since for any $\vf \in \calf C_b^2(H)$, $t\geq 0$, $x\in H$,
\begin{align*}
   \int_H P_t \vf(y)\nu(dy)
   = \lim\limits_{s\to \9} P_s (P_t \vf)(x)
   = \lim\limits_{s\to \9} P_{s+t} \vf(x)
   = \int_H  \vf(y) \nu(dy),
\end{align*}
we have
\begin{align} \label{N0-infi-rde}
   \int_H N_0 \vf(y) \nu(dy)
   = \frac{d}{dt} (\int_H P_t \vf(y) \nu(dy)) |_{t=0}  =0,\ \ \forall \vf \in \calf C_b^2(H),
\end{align}
which implies $(H1)$ $(ii)$.
Moreover, the results of \cite[Section 4.6]{DP04} presented above show that
$(N_0, \mathcal{E}_\Delta(H))$ is essentially m-dissipative,
and so is $(N_0\vf, \calf C_b^2(H))$,
thereby yielding $(H1)$ $(iii)$.
Hence, Hypothesis $(H1)$ is fulfilled.

Concerning Hypothesis $(H1)'$ we have
\begin{proposition} \label{Prop-WP-M-rde}
Assume $(B)$.
Then, Hypothesis $(H1)'$ is satisfied,
i.e.,
the martingale problem for \eqref{equa-rd-X} is well posed in the case $u\equiv 0$.
\end{proposition}

{\bf Proof.}
Set $\ol{H} := L^{2d}(\calo)$.
We have $\nu(\ol{H}) =1$.
For each $x\in \ol{H}$,
by Theorem 4.8 of \cite{DP04},
there exists a unique $(\mathscr{F}_t)$-adapted
process $X(\cdot, x)$,
such that $X(t, x) \in \ol H$ for all $t\geq 0$,
$X\in C([0,T]; L^2(\Omega; H))$,
$\bbe \|X(t,x)\|^{2d}_{L^{2d}(\mathcal{O})}
\leq C_{m,p,T} (1+ \|x\|^{2d}_{L^{2d}})$ for any $m\geq 1$,
and
$X$  solves  \eqref{equa-rd-X} in the mild sense,
i.e.,
for each $t\in [0,T]$,
\begin{align} \label{equa-rde-mild}
   X(t, x) = e^{t\D} x + \int_0^t e^{(t-s) \Delta} F(X(s, x)) ds + W_{\D} (t),\ \ \bbp-a.s.,
\end{align}
where $F(X(s, x)) = -p(X(s, x))$,
\begin{align*}
   W_{\D}(t) = \int_0^t e^{(t-s)\D} (-\D)^{-\frac \g 2} d W(s)
             = \sum\limits_{k=1}^\9 \int_0^t e^{(t-s)\D} (-\D)^{-\frac \g 2} e_k d\beta_k(s),\ t\geq 0,
\end{align*}
and $\{e_k\}$ is the eigenbasis of $-\D$,
i.e., $-\D e_k = \lbb_k e_k$, $k\geq 1$.
Moreover,
$X$ is a Feller process (see \cite[Proposition 4.9]{DP04}).

Note that,
by the integrabilities of $X$ above,
we have that $\bbp$-a.s. $F(X(t, x)) \in H$
and $ t\mapsto \int_0^t F(X(s, x)) ds $ is continuous in $H$.
Taking into account the continuity of $ W_{\D}$ in $H$
implied by \cite[Proposition 4.3]{DP04},
we can take a $\bbp$-version of the process $X$ (still denoted by $X$),
such that $X \in C([0,T]; H)$, $\bbp$-a.s.,
and $X$ satisfies \eqref{equa-rde-mild} for all $t\in [0,T]$ outside a common $\bbp$-null set.
Below we consider this $\bbp$-version process $X$.

Next,
let $x_k := \<x, e_k\>$, $x\in H$, $k \geq 1$.
We claim that $\bbp$-a.s.  for each $k\geq 1$ and for all $t\in [0,T]$,
\begin{align} \label{equa-Xk-rde}
   X_k(t, x) = e^{-\lbb_k t} x_k
            + \int_0^t e^{-\lbb_k (t-s)} (F(X(s, x)))_k ds
            + \int_0^t e^{-\lbb_k (t-s)} \lbb_k^{-\frac \g 2} d\beta_k(s).
\end{align}

For this purpose,
we first infer from \eqref{equa-rde-mild} that $\bbp$-a.s.
for each $k\geq 1$
\begin{align} \label{equa-Xk-rde-mild}
   X_k(t, x) = e^{-\lbb_k t} x_k
           + \l<\int_0^t e^{(t-s)\Delta} F(X(s, x)) ds, e_k\r>
           + \l< W_\D(t), e_k\r>,\ t\in [0,T].
\end{align}
Since $e^{(t-\cdot)\Delta} F(X(\cdot, x))$ is Bochner integrable on $H$
and $y\mapsto \<y, e_k\>$ is a linear bounded operator on $H$,
we get
\begin{align}  \label{F-rde}
   \l<\int_0^t e^{(t-s)\Delta} F(X(s, x)) ds, e_k\r>
   =& \int_0^t   \l<e^{(t-s)\Delta} F(X(s, x)), e_k\r>   ds \nonumber \\
   =& \int_0^t  e^{-\lbb_k(t-s)} (F(X(s, x)))_k   ds
\end{align}
Moreover, since
\begin{align}
   \bbe \sum\limits_{j=1}^\9
   |\l<\int_0^t e^{(t-s)\D} (-\D)^{-\frac \g 2} e_j d \beta_j(s), e_k \r>|^2 <\9,
\end{align}
using Fubini's theorem to exchange the integration with sum
we get
\begin{align} \label{W-rde}
   \l<\sum\limits_{j=1}^\9  \int_0^t e^{(t-s)\D} (-\D)^{-\frac \g 2} e_j d\beta_j(s), e_k \r>
   =& \sum\limits_{j=1}^\9  \l< \int_0^t e^{(t-s)\D} (-\D)^{-\frac \g 2} e_j d\beta_j(s), e_k \r>   \nonumber \\
   =&  \int_0^t   e^{-\lbb_k (t-s)} \lbb_k^{-\frac \g 2}   d\beta_k(s).
\end{align}
Thus, plugging \eqref{F-rde} and \eqref{W-rde} into \eqref{equa-Xk-rde-mild}
we obtain \eqref{equa-Xk-rde}, as claimed.

Hence, we infer from \eqref{equa-Xk-rde}  that $\bbp$-a.s.
\begin{align} \label{equa-Xk-rde-weak}
   d X_k(t, x) = - \lbb_k X_k(t, x) dt
              +  (F(X(t, x)))_k dt
              + \lbb_k^{-\frac \g 2} d\beta_k(t)
\end{align}
with $X_k(0, x) = x_k$.
Since for each $\vf\in \calf  C_b^2(H)$,
there exists $\phi \in C_b^2(\bbr^n)$ such that
$\vf(x) = \phi(\<x,e_1\>,\cdots, \<x,e_n\>)$ for some $n \in \mathbb{N}$,
using It\^o's formula we obtain that,
if $X^n:= (X_1,\cdots, X_n)$,
\begin{align*}
   d\vf(X(t, x))
   =& \sum\limits_{k=1}^n (-\lbb_k  X_k + (F(X(t, x)))_k) \p_k \phi(X^n(t, x)) dt \\
    & + \frac 12 \sum\limits_{k=1}^n \lbb_k^{-\g} \p_{kk}\phi(X^n(t, x)) dt
      + \sum\limits_{k=1}^n \lbb^{-\frac \g 2}_k \p_k \phi(X^n(t, x)) d \beta_k(t) \\
   =& N_0 \vf(X(t, x)) dt + \sum\limits_{k=1}^n \lbb^{-\frac \g 2}_k \p_k \phi(X^n(t, x)) d \beta_k(t),
\end{align*}
This yields that $\vf (X(t, x)) - \int_0^t N_0 \vf(X(s, x)) ds $ is an $(\mathscr{F}_t)$-martingale
under $\bbp$,
so the property $(ii)$ in Definition \ref{Def-equacontr} is fulfilled.

Therefore,
let $\wt \Omega:= C([0,T];H)$,
$\wt{\mathscr{F}}:= \sigma(\mathscr{F}^{X(\cdot,x)}, x\in \ol{H})$
and $\wt{\mathscr{F}}_t:= \sigma(\mathscr{F}_t^{X(\cdot,x)}, x\in \ol{H})$,
$0\leq t\leq T$,
where $\mathscr{F}^{X(\cdot,x)}$ and $\mathscr{F}_t^{X(\cdot,x)}$ denote the
image $\sigma$-algebras under $X(\cdot,x)$ of $\mathscr{F}$ and $\mathscr{F}_t$, respectively.
Set $\pi_t(\omega): =\omega(t)$
and $\wt{\bbp}_x:= \bbp \circ X(\cdot, x)^{-1}$,
$\omega\in \wt \Omega$, $0\leq t\leq T$, $x\in \ol{H}$.
Then,
$(\wt{\Omega}, \wt{\mathscr{F}}, (\wt{\mathscr{F}}_t)_{t\geq 0},
(\pi_t)_{t\geq 0}, (\wt\bbp_x)_{x\in \ol{H}})$
solves the martingale problem for \eqref{equa-rd-X}.
Taking into account Remark \ref{Rem-A-bdd} $(ii)$
we finish the proof of Proposition \ref{Prop-WP-M-rde}.
\hfill $\square$

In the case where $\g=0$,
Hypothesis $(H3)$ holds in certain situations.
Actually,
\cite[Theorem 4.26]{DP04} yields that
$D$ is closable from $L^2(H,\nu)$ to $L^2(H; H, \nu)$,
and it also follows from \cite[Theorem 4.34]{DP04} that
the invariant measure $\nu$ has the density $\rho = \frac{d\nu}{d\mu}$
with respect to the Gaussian
measure $\mu$
with mean zero and covariance operator $-\frac 12 A^{-1}$.
If, in addition, for some $\ve \in (0,1)$,
\begin{align} \label{den-rd}
   \int_H |D \log \rho|_H^{2+\ve} d\nu <\9,
\end{align}
then,
by \cite[Theorem 4.35]{DP04},
$W^{1,2}(H,\nu)$ is compactly embedded into $L^2(H,\nu)$,
and so Hypothesis $(H3)$ is satisfied.

In conclusion, we have from Theorem \ref{Thm-Contr-SDE} that
\begin{theorem} \label{Thm-Contr-SRDE}
Consider the controlled stochastic reaction-diffusion equation \eqref{equa-rd-X}.
Assume   $(B)$.
Assume also $(H2)$ or $(H3)$.
Then, for any $g\in D(N_2)$,
there exists an optimal control $u_*$ for the optimal control problem below
\begin{align*}
  {\rm Min} \bigg\{\int_0^T \int_H &\bbe_{\bbp_x} g(X^u(t)) \nu(dx) dt;\ u\in \calu_{ad},\
                    \ \bbp_x\circ (X^u)^{-1}\ solves\  \nonumber \\
              & the\ martingale\ problem\ for\ \eqref{equa-rd-X}\ for\ \nu-a.e.\ x\in H \bigg\}.
\end{align*}
In particular,
in the case where $\g=0$ and that Assumption $(B)$ and \eqref{den-rd} hold,
we can take $B=Id$.
\end{theorem}

\subsection{Stochastic porous media equations} \label{Subsec-SPME}

In this subsection,
we are concerned with the optimal control problems for stochastic porous media equations.
Precisely,
we consider the controlled stochastic low diffusion equation as in \cite{BBDR06}
\begin{align} \label{equa-X-PME-Low}
   d X(t) &= \Delta (\Psi(X(t))) dt + Q^\frac 12 B u(X(t)) dt + Q^\frac 12 dW(t), \\
   X(0) &= x\in H.   \nonumber
\end{align}
Here $H = H^{-1}(\mathcal{O})$, which is the dual space of $H^1_0(\calo)$
equipped with the inner product
$\<x,y\> := \int_\mathcal{O} ((-\Delta)^{-1} x)(\xi) y(\xi) d\xi$,
$\calo \subseteq \bbr^d$ is a bounded open set
with Dirichlet boundary conditions for the Laplacian $\Delta$,
$B$ and $Q$ are as in \eqref{equa-X},
and $\Psi$ is a dissipative nonlinearity.

In this case,
$A(x) = \Delta (\Psi(x))$
and $D(A) = \{x\in L^2(\calo), \Psi(x) \in H_0^1(\calo)\}$.

Following \cite{BBDR06}, we assume
\begin{enumerate}
  \item[(C1)] There exist $q_k \in [0,\9)$, $k\in \mathbb{N}$,
              such that for the eigenbasis $\{e_k\}$ of $\Delta$ in $H$,
               $Qe_k = q_k e_k$, $k \in \mathbb{N}$.
  \item[(C2)] $\sum_{k=1}^\9 \sup_{\xi\in D}|e_k(\xi)|^2 q_k <\9$.
  \item[(C3)] $\Psi \in C^1(\bbr)$, $\Psi(0)=0$,
              and there exist $r\in (1,\9)$ and   $\kappa_0, \kappa_1, C_1>0$ such that
              \begin{align*}
                    \kappa_0 |s|^{r-1} \leq \Psi'(s) \leq \kappa_1 |s|^{r-1} + C_1, \ \ \forall s\in \bbr.
              \end{align*}
\end{enumerate}

It is known (\cite[Proposition 3.1]{BBDR06}) that,
under Assumption $(C3)$,
$A$ is m-dissipative on $H$.

The corresponding Kolmogorov operator is formally given by
\begin{align*}
  N_0 \vf (x):= \frac 12 \sum\limits_{k=1}^\9 q_k \partial_{kk} \vf(e_k, e_k)
                +  \<\Delta \Psi(x), D\vf(x)\>, \ \ x\in H, \ \vf\in \calf C_b^2(H).
\end{align*}
As in \cite{BBDR06},
let $\mathcal{M}$ be the set of infinitesimally excessive measures $\nu$,
which are infinitesimally invariant measures for $N_0$
satisfying  $(H1)$ $(ii)$,
\begin{align} \label{naPsi}
   \int_H \int_D |\na(\Psi(x))(\xi)|^2 d\xi \nu(dx) <\9,
\end{align}
and for some $\lbb_\nu \in (0,\9)$
\begin{align} \label{N0-lbb-lde}
   \int_H N_0\vf(x) \nu(dx) \leq \lbb_\nu \int_H \vf \nu(dx),
    \ \  \forall \vf \in \calf C_b^2(H)\ with\ \vf\geq 0,\ \nu-a.e..
\end{align}

We see that
$(H1)$  $(i)$ is satisfied for each $\nu \in \mathcal{M}$.
Actually, by Poincar\'{e}'s inequality and \eqref{naPsi},
\begin{align*}
   \int_H \|\Psi(x)\|_{L^2(\calo)}^2 \nu(dx)
   \leq C\int_H \|\na \Psi(x)\|_{L^2(\calo)}^2 \nu(dx) <\9,
\end{align*}
which along with Assumption $(C3)$ above yields that
\begin{align}
   \int_H |x|_H^{2r} \nu(dx)
   \leq \int_H \|x\|_{L^2(\calo)}^{2r} \nu(dx)
   \leq C (1+\int_H \|\Psi(x)\|_{L^2(\calo)}^2 \nu(dx)) <\9.
\end{align}
Moreover,
by \eqref{naPsi},
\begin{align}
   \int_H |\Delta \Psi(x)|_H^2 \nu(dx)
   = \int_H \|\na \Psi(x)\|_{L^2(\calo)}^2 \nu(dx) <\9.
\end{align}
Hence, Hypothesis $(H1)$ $(i)$ follows.

Moreover,
under  Assumptions $(C1)$-$(C3)$,
it follows from \cite[Theorem 4.1]{BBDR06} that
$(N_0, C_b^2(H))$ is essentially m-dissipative on $L^2(H,\nu)$
for each $\nu \in \mathcal{M}$,
which implies $(H1)$ $(iii)$,
and so Hypothesis $(H1)$ holds.

Furthermore,
if in addition $r \geq 2$,
\cite[Theorems 5.1]{BBDR06} yields that
the martingale problem for \eqref{equa-X-PME-Low}
has a solution in the case $u\equiv 0$.
Then, taking into account Remark \ref{Rem-A-bdd} $(ii)$ on uniqueness
we infer that Hypothesis $(H1)'$ holds.

Therefore,
in view of  Theorems \ref{Thm-Contr} and \ref{Thm-Contr-SDE}, we obtain
\begin{theorem}
Consider the controlled stochastic low diffusion equation \eqref{equa-X-PME-Low}.
Assume Hypotheses $(C1)$, $(C2)$ and $(C3)$.
Assume additionally $(H2)$.
Then, for any $g\in D(N_2)$,
there exists an optimal control $u_*$ for the optimal control problem $(P^*)$.

Moreover,
if in addition $r\geq 2$,
we have the optimal controllers for the problem below
\begin{align*}
  {\rm Min} \bigg\{\int_0^T \int_H & \bbe_{\bbp_x}g(X^u(t)) \nu(dx) dt;\ u\in \calu_{ad},\
                     \ \bbp_x\circ (X^u)^{-1}\ solves\ \nonumber \\
              & the\ martingale\ problem\ for\ \eqref{equa-X-PME-Low}\ for\ \nu-a.e.\ x\in H \bigg\}.
\end{align*}
\end{theorem}

{\bf Acknowledgement}
The first and third authors would like to thank the warm hospitality
of Bielefeld University in 2018
where this work was initiated.
The first two authors also thank  the hospitality at Shanghai Jiao Tong University
in 2019 where most part of this work was done.
The third author is supported by NSFC (No. 11871337).
Financial support by the DFG through CRC 1283 is also gratefully acknowledged.


\begin{thebibliography}{99}


\bibitem{B19}
V. Barbu,
Optimal feedback controllers for the stochastic reflection problem in $\bbr^d$,
preprint.


\bibitem{B10}
V. Barbu,
Nonlinear differential equations of monotone types in Banach spaces.
{\it Springer Monographs in Mathematics}.
Springer, New York, 2010. x+272 pp.


\bibitem{BBDR06}
V. Barbu, V.I. Bogachev, G. Da Prato, M. R\"ockner,
Weak solutions to the stochastic porous media equation via Kolmogorov equations: the degenerate case.
{\it J. Funct. Anal}. {\bf 237} (2006), no. 1, 54--75.

\bibitem{BDR16}
V. Barbu, G. Da Prato, M. R\"ockner,
Stochastic porous media equations. {\it Lecture Notes in Mathematics}, {\bf 2163}. Springer, 2016. ix+202 pp.

\bibitem{BBR06}
L. Beznea, N. Boboc, M. R\"ockner,
Markov processes associated with $L^p$-resolvents and applications to stochastic differential equations on Hilbert space. {\it J. Evol. Equ.} {\bf 6} (2006), no. 4, 745--772.

\bibitem{BCR19}
L. Beznea, I. Cimpean, M. R\"ockner,
Treatment of noit allowed starting points for solutions
to singular SDEs on Hilbert spaces: a general approach,
arXiv: 1904.01607v1.



\bibitem{CG95}
A. Chojnowska-Michalik, B. Goldys,
Existence, uniqueness and invariant measures for stochastic semilinear equations on Hilbert spaces.
{\it Probab. Theory Related Fields} {\bf 102} (1995), no. 3, 331--356.

\bibitem{CIL92}
M.G. Crandall, H. Ishii, P.L. Lions,
User's guide to viscosity solutions of second order partial differential equations.
{\it Bull. Amer. Math. Soc.}  {\bf 27} (1992), no. 1, 1--67.

\bibitem{DP04}
G. Da Prato, Kolmogorov equations for stochastic PDEs.
{\it Advanced Courses in Mathematics}. CRM Barcelona.
Birkh\"auser Verlag, Basel, 2004. x+182 pp.


\bibitem{DDG02}
G. Da Prato, A. Debussche, B. Goldys,
Some properties of invariant measures of non symmetric dissipative stochastic systems.
{\it Probab. Theory Related Fields} {\bf 123} (2002), no. 3, 355--380.

\bibitem{DFPR13}
G. Da Prato, F. Flandoli, E. Priola, M. R\"ockner,
Strong uniqueness for stochastic evolution equations in Hilbert spaces perturbed by a bounded measurable drift.
{\it Ann. Probab}. {\bf 41} (2013), no. 5, 3306--3344.

\bibitem{DFPR15}
G. Da Prato, F. Flandoli, E. Priola, M. R\"ockner,
Strong uniqueness for stochastic evolution equations with unbounded measurable drift term.
{\it J. Theoret. Probab}. {\bf 28} (2015), no. 4, 1571--1600.

\bibitem{DFRV16}
G. Da Prato, F. Flandoli, M. R\"ockner, A.Yu Veretennikov,
Strong uniqueness for SDEs in Hilbert spaces with nonregular drift.
{\it Ann. Probab}. {\bf 44} (2016), no. 3, 1985--2023.

\bibitem{DL14}
G. Da Prato, A. Lunardi,
Sobolev regularity for a class of second order elliptic PDE's in infinite dimension.
{\it Ann. Probab}. {\bf 42} (2014), no. 5, 2113--2160.

\bibitem{DR02}
G. Da Prato, M. R\"ockner,
Singular dissipative stochastic equations in Hilbert spaces.
{\it Probab. Theory Related Fields}. {\bf 124} (2002), no. 2, 261--303.

\bibitem{DRW09}
G. Da Prato, M. R\"ockner, F.Y. Wang,
Singular stochastic equations on Hilbert spaces:
Harnack inequalities for their transition semigroups. {\it J. Funct. Anal.} {\bf 257} (2009), no. 4, 992--1017.


\bibitem{FGS17}
G. Fabbri, F. Gozzi, A. \'{S}wiech,
Stochastic optimal control in infinite dimension.
Dynamic programming and HJB equations. With a contribution by Marco Fuhrman and Gianmario Tessitore.
Probability Theory and Stochastic Modelling, 82. Springer, Cham, 2017. xxiii+916 pp.


\bibitem{FT02}
M. Fuhrman, G. Tessitore,
Nonlinear Kolmogorov equations in infinite dimensional spaces:
the backward stochastic differential equations approach and applications to optimal control.
{\it Ann. Probab}. {\bf 30} (2002), no. 3, 1397--1465.


\bibitem{GW02}
F.Z. Gong, F.Y. Wang,
Functional inequalities for uniformly integrable semigroups and application to essential spectrums. {\it Forum Math}. {\bf 14} (2002), no. 2, 293--313.


\bibitem{L69}
J.L. Lions,
Quelques m\'{e}thodes de r\'{e}solution des probl\`{e}mes aux limites non lin\'{e}aires.
(French) Dunod; Gauthier-Villars, Paris 1969 xx+554 pp.

\bibitem{MR92}
Z.M. Ma, M. R\"ockner,
Introduction to the theory of (nonsymmetric) Dirichlet forms. {\it Universitext}. Springer-Verlag, Berlin, 1992. vi+209 pp.

\bibitem{S99}
W. Stannat, The theory of generalized Dirichlet forms and its applications in analysis and stochastics.
{\it Mem. Amer. Math. Soc}. {\bf 142} (1999), no. 678, viii+101 pp.

\bibitem{T03}
G. Trutnau,
On a class of non-symmetric diffusions containing fully nonsymmetric distorted Brownian motions.
{\it Forum Math}. {\bf 15} (2003), no. 3, 409--437.

\bibitem{W00}
F.Y. Wang, Functional inequalities, semigroup properties and spectrum estimates.
{\it Infin. Dimens. Anal. Quantum Probab. Relat. Top}. {\bf 3} (2000), no. 2, 263--295.


\bibitem{W17}
F.Y. Wang, Integrability conditions for SDEs and semilinear SPDEs.
{\it Ann. Probab}. {\bf 45} (2017), no. 5, 3223--3265.




\end{thebibliography}
\end{document}